# Symmetry-Breaking Global Bifurcation in a Surface Continuum Phase-Field Model for Lipid Bilayer Vesicles

**Timothy J. Healey[a*], Sanjay Dharmavaram[b]**

[a]*Department of Mathematics, Cornell University, Ithaca, NY*
[b]*Department of Mechanical and Aerospace Engineering, UCLA, Los Angleles, CA*

---

## Abstract

We study a model for lipid-bilayer membrane vesicles exhibiting phase separation, incorporating a phase field together with membrane fluidity and bending elasticity. We prove the existence of a plethora of equilibria in the large, corresponding to symmetry-breaking solutions of the Euler-Lagrange equations, via global bifurcation from the spherical state. To the best of our knowledge, this constitutes the first rigorous existence results for this class of problems. We overcome several difficulties in carrying this out. Due to inherent surface fluidity combined with finite curvature elasticity, neither the Eulerian (spatial) nor the Lagrangian (material) description of the model lends itself well to analysis. This is resolved via a singularity-free radial-map description, which effectively eliminates the grossly underdetermined in-plane deformation. The resulting governing equations comprise a quasi-linear elliptic system with nonlinear constraints. We show the equivalence of the problem to that of finding the zeros of compact vector field. The latter is not routine. Using spectral and a-priori estimates together with Fredholm properties, we demonstrate that the principal part of the quasi-linear mapping defines an operator with compact resolvent. We then combine well known group-theoretic ideas for symmetry breaking with global bifurcation theory to obtain our results.

*Mathematics Subject Classification* (2010): 35B32 · 35B38 · 35B45 · 35G99 · 35J60 · 35J99 · 35Q99 · 37G40 · 58J55 · 74K25

---

## 1. Introduction

The present work is inspired to a large extent by the striking images presented in [4] of patterned configurations of two-phase fluid domains in man-made lipid-bilayer vesicles. The importance of understanding their mechanical behavior, due to their close relationship with bio-molecular systems and cell function, is well known and highlighted in [4]. Here we analyze a continuum phase-field model for such structures, incorporating curvature elasticity, membrane fluidity and a Cahn-Hilliard type phase-field energy, the latter approximating a sharp-interface theory, e.g., cf. [2], [3], [12], [24], [38]. As in [12], our model features phase-field-dependent bending moduli, which we show enable the identification of new phenomena in the local, asymptotic bifurcation behavior. Specialized to constant bending moduli, the model is the same as that considered in [38], and in the absence of the phase field, the energy reduces to the well-known model of Helfrich [18]. Each of the above-mentioned works provides an informative scientific introduction to the subject (which we do not attempt replicate here) along with numerical results and/or formal analyses of axisymmetric configurations. Here we prove the existence of a plethora of equilibria, corresponding to symmetry-breaking solutions of the Euler-Lagrange equations, via global bifurcation from the spherical state. To the best of our knowledge, this constitutes the first rigorous existence results for this class of problems. Several difficulties arise in this endeavor, as enumerated in the outline that follows.

---

[*]Corresponding author.
  *E-mail address:* healey@math.cornell.edu

In Section 2 we formally derive the Euler-Lagrange equilibrium equations via the first-variation condition for the total potential energy in the presence of constraints. This is a rather elaborate calculation, due in part to the fact that the energy density is naturally prescribed per unit current surface area. Nonetheless, it constitutes a systematic and reliable procedure for obtaining the correct governing equations. Much previous work along these lines for similar models has been carried out elsewhere, e.g., [7],[8], [9], [37[, of which we take full advantage here. In this way we obtain the intrinsic Eulerian description of the equilibrium equations with respect to the deformed surface, in consonance with the membrane fluidity. We impose two constraints: one enforcing area preservation and the other prescribing the average value of the phase field. The latter serves as our bifurcation parameter. Throughout we presume a fixed, nonnegative internal pressure; our philosophy here is to study inflationary phase transitions, as opposed to the more traditional shell buckling induced by compressive external pressure.

Of course the deformed surface is not known a-priori, ostensibly requiring a Lagrangian description for the purpose of analysis. However, due to the aforementioned fluidity, a total Lagrangian formulation is necessarily invariant under the set of all automorphisms of the undeformed reference surface into itself. Such gross non-isolation of solutions lends itself well to neither analysis nor numerics. In Section 3 we employ a singularity-free radial-graph description, previously used in the study of surfaces with prescribed mean curvature [39]. This effectively eliminates the underdetermined in-plane deformation. The resulting equations, written with respect the fixed unit sphere, $S^2$, now appear appreciably more complicated, but nonetheless comprise a quasi-linear, uniformly elliptic system for the unknown phase field and the surface placement field, in the presence of nonlinear constraints.

In Section 4 we demonstrate the equivalence of our problem to an abstract-operator equation for the zeros of a compact vector field. This construction is not routine. The principal part of the quasi-linear system is second order in the phase field and fourth order in the placement field. We use spectral and a-priori estimates together with the stability of the Fredholm index to demonstrate that this defines an operator with compact resolvent [26]. With this in hand, we obtain the equivalent formulation via the inversion of a spectral shift of the principal operator. This sets the stage for global bifurcation methods via the Leray-Schauder degree [35].

We consider the rigorous linearization of the system about the trivial family of spherically symmetric solutions in Section 5. In particular, we determine the basic necessary condition for bifurcation, while characterizing the high-dimensional null space of the linear operator at potential bifurcation points. The latter is a direct consequence of spherical symmetry. Here already we see the role of the bending moduli functions within the so-called spinodal region: Only their non-constancy allows for the generic participation of the displacement variable in the null vectors. Said differently, if the bedning moduli are constant, only the phase field is present in the null vectors asymptotically associated with eventual solutions. All calculations in Section 5 are first conveniently carried out for the pde-based formulation of Section 2. We then demonstrate the equivalence of those results for the abstract formulation of Section 4, the latter of which is needed for the purposes of global bifurcation methods.

In Section 6 we first establish the $O(3)$ equivariance of the field equations. Next we use well-known symmetry ideas to obtain reduced problems amenable to global bifurcation analysis. We exploit the group-theoretic classification in [14], identifying those subgroups that lead to generic bifurcation problem having a one-dimensional null space. Here we work in fixed-point spaces of the Banach space, enabling global conclusions. With this in hand, we demonstrate, for physically reasonable parameter regimes, that *all* of the subgroups catalogued in [14, Thm.9.9] lead to global solution branches, each characterized by precise symmetry. We illustrate the latter in several examples for specific subgroups and mode numbers. In each case we give the explicit linear combination of null vectors (bifurcation direction) asymptotically associated with the local behavior of the solution branch near the bifurcation point, together with a graphical illustration of the associated eigenfunction. These examples are chosen for their apparent relevance to the experimental images shown in [4].

In Section 7 we make a few final remarks. Among other things, we discuss both the applicability of our analysis to other related models, and the attractiveness of our formulation to numerical implementation. We also give a simple argument showing that even in the case of constant bending moduli, there can be no solutions to our problem characterized by a *nontrivial* (non-constant) phase field on the undeformed sphere.

**Notation**

We employ coordinate-free notation as much as possible throughout this work, similar to that introduced in [16], which we now summarize for the convenience of the reader. Let $\Sigma \subset \mathbb{R}^3$ denote a smooth, closed surface, and let $\mathcal{T}_{\mathbf{x}} \subset \mathbb{E}^3$ denote the two-dimensional tangent space at $\mathbf{x} \in \Sigma$, where $\mathbb{E}^3$ is the translation (vector) space associated with point space $\mathbb{R}^3$. Since $\Sigma$ is presumed smooth, for any $\mathbf{x} \in \Sigma$, there is a smooth (locally bijective) map $\mathcal{X}_x : \mathcal{N}_x \to \Sigma$, with $\mathcal{N}_x \subset \mathcal{T}_x$ open, such that $\mathcal{X}_x(\mathbf{0}) = \mathbf{x}$. For a smooth scalar field $\psi(\cdot)$ defined on $\Sigma$, the *surface gradient* at $\mathbf{x} \in \Sigma$, denoted $\nabla_\Sigma \psi(\mathbf{x}) \in \mathcal{T}_\mathbf{x}$, is given by

$$\nabla_\Sigma \psi(\mathbf{x}) := \nabla(\psi \circ \mathcal{X}_x)(\mathbf{0}). \tag{1.1}$$

For a smooth vector field $\mathbf{v} : \Sigma \to \mathbb{E}^3$, the surface gradient at $\mathbf{x} \in \Sigma$, denoted $\nabla_\Sigma \mathbf{v}(\mathbf{x}) \in L(\mathcal{T}_\mathbf{x}, \mathbb{E}^3)$, is given analogously to (1.1). Let $\mathbf{P}_\mathbf{x}$ denote the orthogonal projection onto $\mathcal{T}_x$, viz.,

$$\mathbf{P}_\mathbf{x} = \mathbf{I} - \mathbf{n}(\mathbf{x}) \otimes \mathbf{n}(\mathbf{x}), \tag{1.2}$$

where $\mathbf{I}$ denotes the identity on $\mathbb{E}^3$ and $\mathbf{n}(\cdot)$ is a smooth unit normal field on $\Sigma$, i.e., $\mathbf{n}(\mathbf{x}) \in \mathcal{T}_x^\perp$ for each $\mathbf{x} \in \Sigma$ and $|\mathbf{n}| \equiv 1$. The *tangential surface gradient* of $\mathbf{v}(\cdot)$ at $\mathbf{x} \in \Sigma$ is then defined by

$$D_\Sigma \mathbf{v}(\mathbf{x}) := \mathbf{P}_\mathbf{x} \nabla_\Sigma \mathbf{v}(\mathbf{x}), \tag{1.3}$$

where $D_\Sigma \mathbf{v}(\mathbf{x}) \in L(\mathcal{T}_\mathbf{x}, \mathcal{T}_\mathbf{x})$. The curvature tensor at $\mathbf{x} \in \Sigma$, denoted $\mathbf{L}(\mathbf{x}) \in L(\mathcal{T}_\mathbf{x}, \mathcal{T}_\mathbf{x})$, is given by

$$\mathbf{L}(\mathbf{x}) := -D_\Sigma \mathbf{n}(\mathbf{x}) = -\mathbf{P}_\mathbf{x} \nabla_\Sigma \mathbf{n}(\mathbf{x}) = -\nabla_\Sigma \mathbf{n}(\mathbf{x}),$$

with mean and Gaussian curvature given by $H(\mathbf{x}) := tr\mathbf{L}(\mathbf{x})/2$ and $K(\mathbf{x}) := \det \mathbf{L}(\mathbf{x})$, respectively.

A smooth vector field $\mathbf{t}(\cdot)$ on $\Sigma$ satisfying $\mathbf{t}(\mathbf{x}) \in \mathcal{T}_\mathbf{x}$ for each $\mathbf{x} \in \Sigma$, is called a *tangential vector field*, in which case we define the *surface divergence* of $\mathbf{t}(\cdot)$ at $\mathbf{x} \in \Sigma$ via

$$div_\Sigma \mathbf{t}(\mathbf{x}) := tr[D_\Sigma \mathbf{t}(\mathbf{x})], \tag{1.4}$$

where "tr" denotes the trace. Since $\Sigma$ is closed, the divergence theorem reduces to

$$\int_{\Sigma} div_{\Sigma} \mathbf{t}(\mathbf{x}) ds = 0. \tag{1.5}$$

The surface gradient of a smooth scalar field $\psi(\cdot)$ is naturally tangential, i.e., $\nabla_{\Sigma} \psi(\mathbf{x}) \in \mathcal{T}_{\mathbf{x}}$. By virtue of (2.6), the *second tangential derivative* of $\psi(\mathbf{x})$, denoted $D_{\Sigma}^2 \psi(\mathbf{x}) \in BL(\mathcal{T}_{\mathbf{x}} \times \mathcal{T}_{\mathbf{x}}, \mathcal{T}_{\mathbf{x}})$, is the bilinear operator

$$D_{\Sigma}^2 \psi(\mathbf{x}) := D_{\Sigma}(D_{\Sigma} \psi(\mathbf{x})) = \mathbf{P}_x \nabla_{\Sigma}(\nabla_{\Sigma} \psi(\mathbf{x})), \tag{1.6}$$

where in accordance with (2.6), we have introduced

$$D_{\Sigma} \psi(\mathbf{x}) := \mathbf{P}_x \nabla_{\Sigma} \psi(\mathbf{x}) = \nabla_{\Sigma} \psi(\mathbf{x}), \tag{1.7}$$

i.e., the surface gradient and the tangential surface gradient of a scalar field are one and the same. The *Laplace-Beltrami operator* acting on $\psi(\cdot)$ at $\mathbf{x} \in \Sigma$, is given by

$$\Delta_{\Sigma} \psi(\mathbf{x}) := tr[D_{\Sigma}^2 \psi(\mathbf{x})] = div_{\Sigma} \nabla_{\Sigma} \psi(\mathbf{x}). \tag{1.8}$$

Higher tangential derivatives, evaluated at $\mathbf{x} \in \Sigma$, are multilinear functions of Cartesian products of $\mathcal{T}_{\mathbf{x}}$ into $\mathcal{T}_{\mathbf{x}}$, defined in the obvious way by consecutive composition. In particular,

$$\begin{aligned} D_{\Sigma}^3 \psi(\mathbf{x}) &:= D_{\Sigma}(D_{\Sigma}^2 \psi(\mathbf{x})) = D_{\Sigma}^2(\nabla_{\Sigma} \psi(\mathbf{x})), \\ D_{\Sigma}^4 \psi(\mathbf{x}) &:= D_{\Sigma}(D_{\Sigma}^3 \psi(\mathbf{x})) = D_{\Sigma}^3(\nabla_{\Sigma} \psi(\mathbf{x})). \end{aligned} \tag{1.9}$$

## 2. Intrinsic Eulerian Equilibrium Equations

We begin with the following phase-field elastic-shell *potential energy* functional for a vesicle, written with respect to the *current configuration*, denoted $\Sigma$, and presumed isomorphic to the unit sphere $S^2 \subset \mathbb{R}^3$:

$$\int_{\Sigma} [B(\phi)H^2 + E(\phi)K + \frac{\varepsilon}{2}|\nabla_{\Sigma} \phi|^2 + W(\phi)] \, ds - pV(\Sigma), \tag{2.1}$$

subject to the constraints

$$\begin{aligned} \int_{\Sigma} ds &= 4\pi, \\ \int_{\Sigma} \phi ds &= 4\pi\lambda. \end{aligned} \tag{2.2}$$

Here $H$ and $K$ denote the scalar-valued mean-curvature and Gaussian curvature fields of $\Sigma$, respectively, $\phi$ denotes the scalar *phase field, $B(\cdot), E(\cdot)$* are real-valued *bending moduli* functions, $\varepsilon > 0$ is a small "interfacial" parameter, $p \geq 0$ is the prescribed internal *pressure, $V(\Sigma)$* denotes the volume enclosed by the surface $\Sigma$, $\lambda > 0$ is a "control" parameter, and $W(\cdot)$ is a non-convex potential function. We assume throughout that $W(\cdot)$ is of class $C^3$ and $B(\cdot)$ and $E(\cdot)$ are bounded $C^4$ functions, on $\mathbb{R}$. We further require $B(\cdot)$ to be positive, say, according to

$$B : \mathbb{R} \to [\varepsilon, \infty), \qquad (2.3)$$

and we assume that $W$ is a non-negative, non-convex potential function. Typical graphs of $W, W'$ and $W''$ that we consider here are sketched below in Figure 2.1. In particular, we assume that $W''$ has precisely two zeros, $m_1, m_2$, as shown. The numerical values of $\alpha < m_1 < m_2 < \beta$ play no role in the subsequent analysis. The interval $(m_1, m_2)$ is called the *spinodal region.*

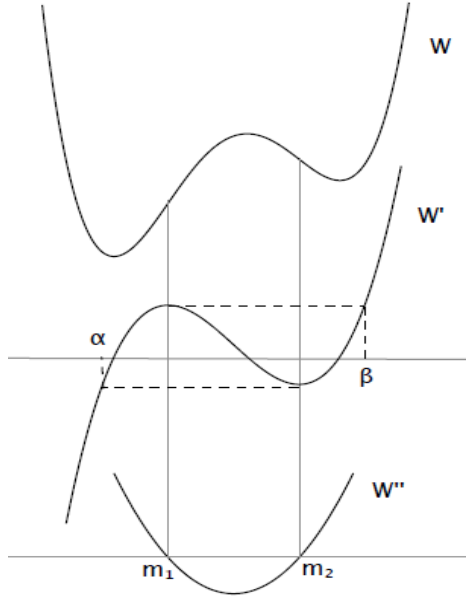

Figure 2.1 Graphs of $W, W'$ and $W''$.

Clearly the first constraint equation $(2.2)_1$ specifies total-area preservation, while $(2.2)_2$ prescribes the average value of the phase field.

Our goal in this section is a formal derivation of the Euler-Lagrange equilibrium equations. This is a rather involved, nontrivial procedure, due in large part to the fact that (2.1) is naturally posed over the current configuration $\Sigma$. Fortunately, much previous work has been carried out along these lines, e.g., [7], [8], [9], [37], of which we take full advantage here. We first introduce a smooth *deformation*

$$\mathbf{f} : S^2 \to \mathbb{R}^3, \Sigma := \mathbf{f}(S^2), \qquad (2.4)$$

denoting points $\mathbf{y} \in \Sigma$ via $\mathbf{y} = \mathbf{f}(\mathbf{x})$ for all $\mathbf{x} \in S^2$. Henceforth we refer to $S^2$ as the *reference configuration*. Next we consider the smooth *variations*

$$\begin{aligned}
\mathbf{f}(\mathbf{x}) &\to \mathbf{f}(\mathbf{x}) + \boldsymbol{\eta}(\mathbf{x}) \Rightarrow \mathbf{y} \to \mathbf{y} + \mathbf{v}(\mathbf{y}) + w(\mathbf{y})\mathbf{n}(\mathbf{y}), \\
\phi(\mathbf{y}) &\to \phi(\mathbf{f}(\mathbf{x}) + \boldsymbol{\eta}(\mathbf{x})) + \psi(\mathbf{y}) = \phi(\mathbf{y}) + \nabla_\Sigma \phi(\mathbf{y}) \cdot \mathbf{v}(\mathbf{y}) + \psi(\mathbf{y}) + ...,
\end{aligned} \qquad (2.5)$$

where $\boldsymbol{\eta} : S^2 \to \mathbb{R}^3$ is the (material) variation of $\mathbf{f}$, $\mathbf{n}$ denotes the outward unit normal field on $\Sigma$, $\mathbf{v} \cdot \mathbf{n} \equiv 0$ on $\Sigma$, and $\psi : \Sigma \to \mathbb{R}$ is the variation of the scalar field $\phi$ on $\Sigma$. Hence, $\mathbf{v}$ and $w\mathbf{n}$ represent the tangential and normal components, respectively, of the variation of the deformed surface in spatial form, and $\psi(\mathbf{y}) + \nabla_\Sigma \phi(\mathbf{y}) \cdot \mathbf{v}(\mathbf{y})$ is the total variation of $\phi(\mathbf{y})$ in spatial form, the second (convected) term accounting for the variation of the deformed surface itself.

We now change variables, expressing (2.1) with respect to the reference sphere, while introducing Lagrange multipliers enforcing the two constraints (2.2):

$$U := \int_{S^2} [B(\phi)H^2 + E(\phi)K + \frac{\varepsilon}{2}|\nabla_\Sigma \phi|^2 + W(\phi)] \, Jds - pV(\Sigma)$$
$$- \gamma \left( \int_{S^2} Jds - 4\pi \right) - \mu \left( \int_{S^2} \phi Jds - 4\pi\lambda \right), \tag{2.6}$$

where $J$ denotes the local area ratio. By a convenient abuse of notation, we employ the same name for the function whether a function of $\mathbf{x}$ or $\mathbf{y}$, e.g., $H(\mathbf{x})$ instead of $H(\mathbf{f}(\mathbf{x}))$, etc. Then the first variation is given formally by

$$\delta U = \int_{S^2} \{2B(\phi)H\delta H + E(\phi)\delta K + \varepsilon \nabla_\Sigma \phi \cdot \delta(\nabla_\Sigma \phi)$$
$$+ [B'(\phi)H^2 + E'(\phi)K + W'(\phi)]\delta\phi\} Jds$$
$$+ \int_{S^2} [B(\phi)H^2 + E(\phi)K + \frac{\varepsilon}{2}|\nabla_\Sigma \phi|^2 + W(\phi)]\delta Jds$$
$$- p\delta V(\Sigma) - \gamma \int_{S^2} \delta Jds - \mu \int_{S^2} (\delta\phi J + \phi\delta J)ds. \tag{2.7}$$

Next we push (2.7) forward to the current configuration and employ the following well known expressions for the variations in spatial form:

$$\delta\phi = \psi + \nabla_\Sigma \phi \cdot \mathbf{v}$$
$$\delta(\nabla_\Sigma \phi) = \nabla_\Sigma(\delta\phi) + [w\mathbf{L} - D_\Sigma \mathbf{v}^T]\nabla_\Sigma \phi,$$
$$\delta H = \nabla_\Sigma H \cdot \mathbf{v} + \Delta_\Sigma w/2 + w(2H^2 - K),$$
$$\delta K = \nabla_\Sigma K \cdot \mathbf{v} - div_\Sigma(Cof\mathbf{L}[\nabla_\Sigma w]) + 2HKw,$$
$$\delta J = J(div_\Sigma \mathbf{v} - 2wH),$$
$$\delta V(\Sigma) = \int_\Sigma wds, \tag{2.8}$$

where $Cof\mathbf{L} \in L(\mathcal{T}_y, \mathcal{T}_y)$ denotes the cofactor tensor corresponding to the curvature tensor $\mathbf{L}$. By virtue of the Cayley-Hamilton theorem, we have

$$Cof\mathbf{L} = 2H\mathbf{I}_y - \mathbf{L}, \tag{2.9}$$

where $\mathbf{I}_y$ denotes the identity on $\mathcal{T}_y$.

Substituting (2.8) into the push forward of (2.7) leads to the first variation condition:

$$\delta U = \int_{\Sigma} \Big\{ [B'(\phi)H^2 + E'(\phi)K + W'(\phi) - \mu]\psi + \varepsilon \nabla_{\Sigma}\phi \cdot \nabla_{\Sigma}\psi$$

$$+ \nabla_{\Sigma}[B(\phi)H^2 + E(\phi)K + W(\phi)] \cdot \mathbf{v}$$

$$+ [B(\phi)H^2 + E(\phi)K + W(\phi)]div_{\Sigma}\mathbf{v}$$

$$+ \varepsilon \Big( (|\nabla_{\Sigma}\phi|^2 / 2)div_{\Sigma}\mathbf{v} + \nabla_{\Sigma}\phi \cdot (\nabla_{\Sigma}[\nabla_{\Sigma}\phi \cdot \mathbf{v}]) - \nabla_{\Sigma}\phi \cdot (D_{\Sigma}\mathbf{v}\nabla_{\Sigma}\phi) \Big)$$

$$- \mu(\nabla_{\Sigma}\phi \cdot \mathbf{v} + \phi div_{\Sigma}\mathbf{v}) - \lambda div_{\Sigma}\mathbf{v}$$

$$+ 2B(\phi)H[\Delta_{\Sigma}w / 2 + (2H^2 - K)w]$$

$$- E(\phi)[div_{\Sigma}([Cof\mathbf{L}]\nabla_{\Sigma}w) - 2HKw]$$

$$+ [\varepsilon(\nabla_{\Sigma}\phi \cdot \mathbf{L}[\nabla_{\Sigma}\phi]) - H|\nabla_{\Sigma}\phi|^2) - p]w$$

$$- 2H[B(\phi)H^2 + E(\phi)K + W(\phi) - \mu\phi - \gamma]w \Big\}ds = 0, \tag{2.10}$$

for all smooth admissible variations $\psi$, $\mathbf{v}$ and $w$. Next we formally integrate by parts and use the divergence theorem (1.5), bringing the entire integrand into the usual form needed to obtain the Euler Lagrange equations. Remarkably this results in the vanishing of the entire tangential term. This is self-evident for all but the term on the fourth line of (2.10), which also vanishes by direct calculation. In any case, the usual localization arguments applied to (2.10) after integrating by parts yields the Euler Lagrange equilibrium equations:

$$-\varepsilon\Delta_{\Sigma}\phi + B'(\phi)H^2 + E'(\phi)K + W'(\phi) = \mu \text{ on } \Sigma, \tag{2.11}$$

$$\Delta_{\Sigma}(B(\phi)H) - 2H\Delta_{\Sigma}(E(\phi)) + (\mathbf{L} \cdot D_{\Sigma}^2(E(\phi))$$

$$+ \varepsilon(\nabla_{\Sigma}\phi \cdot \mathbf{L}[\nabla_{\Sigma}\phi] - H|\nabla_{\Sigma}\phi|^2) + 2B(\phi)H(H^2 - K)$$

$$- 2H[W(\phi) - \gamma - \mu\phi] = p \text{ on } \Sigma, \tag{2.12}$$

subject to the two constraints (2.2). In (2.12) we use the notation $\mathbf{A} \cdot \mathbf{B} := tr(\mathbf{AB}^T)$ for all $\mathbf{A}, \mathbf{B} \in L(\mathcal{T}_y, \mathcal{T}_y)$. Also, in obtaining (2.12) from (2.10) we make use of (2.9) and the fact that $div_{\Sigma}([Cof\mathbf{L}]\mathbf{t}) = Cof\mathbf{L} \cdot D_{\Sigma}\mathbf{t}$ for all tangential fields $\mathbf{t}(\cdot)$ on $\Sigma$, cf. [37].

Equations (2.11), (2.12) elegantly represent the intrinsic *Eulerian* form of the equilibrium conditions, in consonance with the membrane fluidity inherent in the model and given with respect to the unknown surface $\Sigma$. The latter characteristic clearly shows that equations (2.11), (2.12) alone are insufficient. If we introduce a general *vector*-valued mapping in order to describe $\Sigma$, we arrive at the same difficulty associated with a *Lagrangian* description (2.4), where $\mathbf{y} = \mathbf{f}(\mathbf{x})$ is interpreted as the position vector of the material point in $\Sigma$ that occupies position $\mathbf{x} \in S^2$. Let $\mathbf{m} : S^2 \to S^2$ be any smooth automorphism. Then the equilibrium equations (2.11), (2.12) are invariant under $\mathbf{f}(\cdot) \to \mathbf{f} \circ \mathbf{m}(\cdot)$. As a consequence, the tangential component of the deformation $\mathbf{f}$, relative to $S^2$, is parametrized by $\mathbf{m}$ and is thus grossly underdetermined. Of course this kind of indeterminacy also occurs in the Lagrangian description of classical bulk fluids. We conclude that the membrane fluidity is conducive to the Eulerian description, which is at odds with the usual description of the bending elasticity also characterizing our model.

**Remark 2.1** In our recent work [10]**,** we also consider the model (2.1) with (2.2)$_1$ replaced by the local-area constraint $J \equiv 1$, i.e., the sixth term in (2.6) is replaced by

$$\int_{S^2} \gamma_l (J-1) ds,$$

where $\gamma_l$ is now a Lagrange-multiplier *field.* This is similar to the situation in classical fluid mechanics, where the local volume ratio is preserved. The Euler Lagrange equations are essentially the same in this case, the only difference being that the tangential equation does not vanish identically but rather yields the condition $\nabla_{\Sigma} \gamma_l \equiv 0 \Rightarrow \gamma_l = \text{const}$. This strongly hints that the two formulations are "equivalent". However, the distinct constraint equations for the two formulations need to be addressed as well. Existence presumed, we show in [10] that within the equivalence class of (genus-zero) solutions for the global-area formulation, as indicated above in the context of the Lagrangian formulation, there is at least one member that satisfies the local area problem as well. In addition, we show that the formal characterization of the second-variation conditions for the two formulations are equivalent in that sense as well.

## 3. A Singularity-Free Radial-Graph Description

To overcome the above-mentioned shortcomings of the two classical descriptions, we borrow an idea from [39], restricting ourselves to mappings of the form

$$\mathbf{y} = \exp(u(\mathbf{x}))\mathbf{x}, \text{ for all } \mathbf{x} \in S^2, \tag{3.1}$$

where $u : S^2 \to \mathbb{R}$. The non-negative scalar field $\exp(u(\mathbf{x}))$ represents the magnitude of the radial position vector of the deformed surface. While specifying $\Sigma$, we note that $\mathbf{y}$ given by (3.1) does not generally represent the position vector of the material particle occupied by $\mathbf{x}$ in the reference configuration $S^2$. In particular, (3.1) does not specify the tangential component of the displacement. We now proceed to obtain all quantities needed to reformulate our problem, (2.11), (2.12) and (2.2), solely in terms of $\phi$ and $u$ as fields on $S^2$ as follows.

All calculations hereafter are carried out with respect to the reference surface $S^2$. Accordingly, we drop the subscript notation for derivatives, as introduced in Section 1, with subscript $S^2$ now being understood, e.g., $\nabla(\cdot) := \nabla_{S^2}(\cdot)$. We first observe that the outward unit normal field on $S^2$, denoted $\mathbf{n}_o(\cdot)$, is simply the unit-radial translation of the position vector itself, viz., $\mathbf{n}_o(\mathbf{x}) \cong \mathbf{x}$, which we now consistently employ. Let $\mathbf{F}(\mathbf{x}) \in L(T_x, \mathbb{E}^3)$ denote the surface gradient of the right side of (3.1), given by

$$\mathbf{F} = \exp(u)(\mathbf{1}_x + \mathbf{x} \otimes \nabla u), \tag{3.2}$$

where $\mathbf{1}_x$ denotes the identity on $T_x$, the tangent space of $S^2$ at $\mathbf{x}$. Let $\{\mathbf{e}_1, \mathbf{e}_2, \mathbf{x}\}$ denote a right-handed orthonormal frame at $\mathbf{x} \in S^2$. Then

$$\mathbf{a}_\alpha = \mathbf{F}\mathbf{e}_\alpha = \exp(u)(\mathbf{e}_\alpha + u_{,\alpha} \mathbf{x}), \ \alpha = 1, 2, \tag{3.3}$$

define a convected basis, under (3.1), for the tangent plane $\mathcal{T}_y$ on $\Sigma$ at $\mathbf{y}$, where $u_{,\alpha} := \nabla u \cdot \mathbf{e}_\alpha$. From (3.3) we may then compute the outward unit normal, denoted $\mathbf{n}$, to $\Sigma$ at $\mathbf{y}$ as a function of $\mathbf{x} \in S^2$ :

$$\mathbf{n} = (1 + |\nabla u|^2)^{-1/2}(\mathbf{x} - \nabla u). \tag{3.4}$$

The reciprocal basis vectors are given by

$$\begin{aligned}
\mathbf{a}^1 &= \exp(-u)(1 + |\nabla u|^2)^{-1}\{[1 + (u_{,2})^2]\mathbf{e}_1 - u_{,1}\,u_{,2}\,\mathbf{e}_2 + u_{,1}\,\mathbf{x}\}, \\
\mathbf{a}^2 &= \exp(-u)(1 + |\nabla u|^2)^{-1}\{-u_{,1}\,u_{,2}\,\mathbf{e}_1 + [1 + (u_{,1})^2]\mathbf{e}_2 + u_{,2}\,\mathbf{x}\}.
\end{aligned} \tag{3.5}$$

From (3.3)-(3.5) we can directly compute all necessary quantities associated with the deformed surface $\Sigma$. In what follows, Greek indices range from 1 to 2, with repeated indices implying summation:

$$\begin{aligned}
1^{st}\ &\text{fundamental form: } a_{\alpha\beta} = \mathbf{a}_\alpha \cdot \mathbf{a}_\beta; \\
&\mathbf{C} := \mathbf{F}^T\mathbf{F} = a_{\alpha\beta}\mathbf{e}_\alpha \otimes \mathbf{e}_\beta = \exp(2u)(\mathbf{1}_x + \nabla u \otimes \nabla u), \\
&J = \sqrt{a} = (\det[a_{\alpha\beta}])^{1/2} = \exp(2u)(1 + |\nabla u|^2)^{1/2};
\end{aligned} \tag{3.6}$$

$$\begin{aligned}
&a^{\alpha\beta} = \mathbf{a}^\alpha \cdot \mathbf{a}^\beta; \\
&\mathbf{C}^{-1} = a^{\alpha\beta}\mathbf{e}_\alpha \otimes \mathbf{e}_\beta = \exp(-2u)(1 + |\nabla u|^2)^{-1}\mathbf{A}(\nabla u), \\
&\mathbf{A}(\nabla u) := (1 + |\nabla u|^2)\mathbf{1}_x - \nabla u \otimes \nabla u;
\end{aligned} \tag{3.7}$$

$$\begin{aligned}
2^{nd}\ &\text{fundamental form:} \\
&L_{\alpha\beta} = \mathbf{a}_\alpha \cdot \mathbf{L}\mathbf{a}_\beta = \exp(u)(1 + |\nabla u|^2)^{-1/2}(u_{,\alpha\beta} - u_{,\alpha}\,u_\beta - \delta_{\alpha\beta}),
\end{aligned} \tag{3.8}$$

where $u_{,\alpha\beta} := \mathbf{e}_\alpha \cdot [D^2 u]\mathbf{e}_\beta$ and $\delta_{\alpha\beta}$ is the Kronecker delta. From (3.7), (3.8) we then find

$$H = \exp(-u)(1 + |\nabla u|^2)^{-3/2}[\mathbf{A}(\nabla u) \cdot D^2 u - 2(1 + |\nabla u|^2)] / 2; \tag{3.9}$$

$$K = \exp(-2u)(1 + |\nabla u|^2)^{-2}\left(\det D^2 u - \mathbf{A}(\nabla u) \cdot D^2 u + |\nabla u|^2 + 1\right). \tag{3.10}$$

As in Section 2, we conveniently use the same name for the dependent variable as a function of either $\mathbf{x}$ or $\mathbf{y}$; e.g., $H(\mathbf{x})$ given in (3.9), strictly speaking, stands for $H(\exp(u(\mathbf{x}))\mathbf{x})$, etc. It's not hard to see that the symmetric transformation $\mathbf{A}(\nabla u) \in L(T_x, T_x)$ defined in (3.7)$_3$, is positive definite, having eigenvalues $1,\ 1 + |\nabla u|^2$. Hence, (3.9) constitutes a second-order, quasi-linear elliptic equation in $u$, as observed in [39].

In view of (2.11), (2.12) there are a few other quantities that we need to compute. We first observe from (3.2) and (3.3) that the surface gradient is equivalent to

$$\mathbf{F} = \mathbf{a}_\alpha \otimes \mathbf{e}_\alpha, \tag{3.11}$$

and it is convenient to introduce the *restricted inverse* $\mathbf{F}^{-1} \in L(\mathcal{T}_y, T_x)$, defined by

$$\mathbf{F}^{-1} := \mathbf{e}_\alpha \otimes \mathbf{a}^\alpha = \exp(-u)(1 + |\nabla u|^2)^{-1}[\mathbf{A}(\nabla u) + \nabla u \otimes \mathbf{x}], \tag{3.12}$$

where $\mathcal{T}_y$ denotes the tangent plane to $\Sigma$ at $\mathbf{y} = \mathbf{f}(\mathbf{x})$. So defined, observe that $\mathbf{F}^{-1}\mathbf{a}_\alpha = \mathbf{e}_\alpha$, $\alpha = 1, 2$, which is indeed the inverse of (3.3), restricted to tangent planes. From the chain rule for a smooth scalar field $\psi(\cdot)$, we deduce

$$\nabla \psi = \mathbf{F}^T \nabla_\Sigma \psi \Rightarrow \nabla_\Sigma \psi = \mathbf{F}^{-T} \nabla \psi, \tag{3.13}$$

where from (3.11) and (3.12) we have

$$\mathbf{F}^{-T} = \mathbf{a}^\alpha \otimes \mathbf{e}_\alpha = \exp(-u)(1 + |\nabla u|^2)^{-1}[\mathbf{A}(\nabla \mathbf{u}) + \mathbf{x} \otimes \nabla u]. \tag{3.14}$$

Equations (3.6), (3.7), (3.12) and (3.13) yield

$$\begin{aligned}
|\nabla_\Sigma \psi|^2 &= \nabla \psi \cdot (\mathbf{F}^{-1}\mathbf{F}^{-T}\nabla \psi) = \nabla \psi \cdot (\mathbf{C}^{-1}\nabla \psi) \\
&= \nabla \psi \cdot ([a^{\alpha\beta}\mathbf{e}_\alpha \otimes \mathbf{e}_\beta]\nabla \psi),
\end{aligned}$$

and from (3.7) we arrive at

$$|\nabla_\Sigma \psi|^2 = \exp(-2u)(1 + |\nabla u|^2)^{-1}\nabla \psi \cdot ([\mathbf{A}(\nabla u)]\nabla \psi). \tag{3.15}$$

For smooth vector field $\mathbf{v}(\cdot)$, the chain rule gives

$$\nabla \mathbf{v} = [\nabla_\Sigma \mathbf{v}]\mathbf{F} \Rightarrow \nabla_\Sigma \mathbf{v} = [\nabla \mathbf{v}]\mathbf{F}^{-1}. \tag{3.16}$$

In particular, for the vector field $\nabla_\Sigma \psi = \mathbf{F}^{-T}\nabla \psi$, we have

$$\begin{aligned}
\nabla_\Sigma(\nabla_\Sigma \psi) &= \nabla_\Sigma(\mathbf{F}^{-T}\nabla \psi) = \nabla(\mathbf{F}^{-T}\nabla \psi)\mathbf{F}^{-1} \\
&= \mathbf{F}^{-T}\nabla^2\psi\mathbf{F}^{-1} + (\nabla\mathbf{F}^{-T}\nabla \psi)\mathbf{F}^{-1},
\end{aligned} \tag{3.17}$$

where $\nabla^2\psi := \nabla(\nabla\psi)$. By virtue of (3.12), the principal part of (3.17)$_2$ reads

$$\begin{aligned}
\mathbf{F}^{-T}\nabla^2\psi\mathbf{F}^{-1} &= (\mathbf{a}^\alpha \otimes \mathbf{e}_\alpha)\nabla^2\psi(\mathbf{e}_\beta \otimes \mathbf{a}^\beta) \\
&= (\mathbf{e}_\alpha \cdot \nabla^2\psi\mathbf{e}_\beta)\mathbf{a}^\alpha \otimes \mathbf{a}^\beta \\
&= \mathbf{F}^{-T}D^2\psi\mathbf{F}^{-1},
\end{aligned} \tag{3.18}$$

which also shows that $\mathbf{F}^{-T}\nabla^2\psi\mathbf{F}^{-1} = \mathbf{F}^{-T}D^2\psi\mathbf{F}^{-1} \in L(\mathcal{T}_y, \mathcal{T}_y)$. From (1.6)-(1.8) and (3.7), we find that the principal part of the Laplace Beltrami operator is given by

$$\begin{aligned}
\Delta_\Sigma \psi &= tr[\mathbf{P}_y \nabla_\Sigma^2 \psi] + \dots \\
&= tr[\mathbf{F}^{-T}D^2\psi\mathbf{F}^{-1}] + \dots = tr[D^2\psi\mathbf{C}^{-1}] + \dots \\
&= \exp(-2u)(1 + |\nabla u|^2)^{-1}\mathbf{A}(\nabla u) \cdot D^2\psi + \dots
\end{aligned} \tag{3.19}$$

Equations (3.17) and (3.18) show that the lower-order term in (3.19) is given by the expression $tr\left(\mathbf{P}_y[(\nabla\mathbf{F}^{-T}\nabla\psi)\mathbf{F}^{-1}]\right)$, the detailed and lengthy form of which is not needed for our purposes here. We claim that it can be expressed as

$$tr\left(\mathbf{P}_y[(\nabla\mathbf{F}^{-T}\nabla\psi)\mathbf{F}^{-1}]\right) = \exp(-2u)(1+|\nabla u|^2)^{-1}g(D^2u,\nabla u,u)\nabla u \cdot \nabla\psi, \qquad (3.20)$$

where $g:L(T_x,T_x)\times T_x\times\mathbb{R}\to\mathbb{R}$ is smooth. To see this, first note from (1.2) and (3.4) that

$$\mathbf{P}_y = \mathbf{I} - (1+|\nabla u|^2)^{-1}[(\mathbf{x}-\nabla u)\otimes(\mathbf{x}-\nabla u)]. \qquad (3.21)$$

Also, a calculation yields

$$\nabla\mathbf{F}^{-T}\nabla\psi = \exp(-u)\Big\{-\nabla u\otimes\nabla\psi$$
$$+ (1+|\nabla u|^2)^{-1}\big[(\nabla u + 2(1+|\nabla u|^2)^{-1}[D^2u]\nabla u)\otimes\nabla u$$
$$+ D^2u - \mathbf{x}\otimes\nabla u\big](\nabla u \cdot \nabla\psi)\Big\}, \qquad (3.22)$$

where we make use of the fact that $\nabla^2u = D^2u - \mathbf{x}\otimes\nabla u$. Substituting (3.12), (3.21) and (3.22) into the left side of (3.20) gives the form on the right side.

Finally, there are two additional terms in (2.12) that we need to address. First, using (3.7), (3.8), (3.12) and (3.13) we obtain

$$\nabla_\Sigma\phi \cdot \mathbf{L}[\nabla_\Sigma\phi] = \mathbf{F}^{-1}\nabla\phi \cdot L_{\alpha\beta}\mathbf{a}^\alpha\otimes\mathbf{a}^\beta[\mathbf{F}^{-1}\nabla\phi]$$
$$= L_{\alpha\beta}a^{\alpha\gamma}a^{\beta\delta}\phi_{,\gamma}\,\phi_{,\delta}$$
$$= \exp(-3u)(1+|\nabla u|^2)^{-5/2}\nabla\phi\cdot\big([\mathbf{A}(\nabla u)(D^2u-\nabla u\otimes\nabla u-\mathbf{I}_x)\mathbf{A}(\nabla u)]\nabla\phi\big). \qquad (3.23)$$

From (3.18) and (3.19) we find that

$$\mathbf{L}\cdot D_\Sigma^2\psi = \mathbf{L}\cdot[(\mathbf{F}^{-T}D^2\psi + \mathbf{P}_y\nabla\mathbf{F}^{-T}\nabla\psi)\mathbf{F}^{-1}]. \qquad (3.24)$$

The first of these is similar to (3.23) in that

$$\mathbf{L}\cdot(\mathbf{F}^{-T}D^2\psi\mathbf{F}^{-1}) = L_{\alpha\beta}\mathbf{a}^\alpha\otimes\mathbf{a}^\beta\cdot(\mathbf{F}^{-T}D^2\psi\mathbf{F}^{-1})$$
$$= L_{\alpha\beta}a^{\alpha\gamma}a^{\beta\delta}\psi_{,\gamma\delta}$$
$$= \exp(-3u)(1+|\nabla u|^2)^{-1/2}\left(\mathbf{A}(Du)[D^2u-Du\otimes Du-\mathbf{1}_x]\mathbf{A}(Du)\right)\cdot D^2\psi. \qquad (3.25)$$

The lower-order term in (3.24) is more complicated, but as is the case for (3.20), we use (3.12), (3.21) and (3.22) to obtain the representation

$$\mathbf{L}\cdot(\mathbf{P}_y\nabla\mathbf{F}^{-T}\nabla\psi\mathbf{F}^{-1}) = \exp(-2u)(1+|\nabla u|^2)^{-1}\mathbf{r}(D^2u,\nabla u,u)\cdot\nabla\psi, \qquad (3.26)$$

where $\mathbf{r}:L(T_x,T_x)\times T_x\times\mathbb{R}\to T_x$ is smooth and satisfies $\mathbf{r}(\mathbf{0},\mathbf{0},0) = \mathbf{0}$.

We now substitute (3.15), (3.19), (3.20), (3.23), (3.25) and (3.26) into the equilibrium equations (2.11) and (2.12) to obtain

$$-\varepsilon[\mathbf{A}(Du) \cdot D^2\phi + g(D^2u, Du)Du \cdot D\phi]$$
$$+ \exp(2u)(1 + |Du|^2)[B'(\phi)H^2 + E'(\phi)K + W'(\phi) - \mu] = 0 \text{ on } S^2, \tag{3.27}$$

and

$$\mathbf{A}(Du) \cdot \left[ D^2\big(B(\phi)H\big) - 2HD^2\big(E(\phi)\big)\right]$$
$$+ \exp(-u)(1 + |Du|^2)^{1/2}\big(\mathbf{A}(Du)[D^2u - Du \otimes Du - \mathbf{1}_x]\mathbf{A}(Du)\big) \cdot D^2\big(E(\phi)\big)$$
$$+ g(D^2u, Du)Du \cdot \left[ D\big(B(\phi)H\big) - 2HD\big(E(\phi)\big)\right] + \mathbf{r}(D^2u, Du, u) \cdot D\big(E(\phi)\big)$$
$$+ \varepsilon\left\{ \exp(-u)(1 + |Du|^2)^{-3/2} D\phi \cdot \big([\mathbf{A}(Du)[D^2u - Du \otimes Du - \mathbf{1}_x]\mathbf{A}(Du)]D\phi\big)\right.$$
$$\left. - HD\phi \cdot (\mathbf{A}(Du)D\phi)\right\}$$
$$+ \exp(2u)(1 + |Du|^2)\left\{ 2H[B(\phi)(H^2 - K) + \gamma + \mu\phi - W(\phi)] - p\right\} = 0 \text{ on } S^2, \tag{3.28}$$

respectively, where we use observation (1.7). In components, the second-order operator appearing in both (3.27) and (3.28) reads

$$\mathbf{A}(Du) \cdot D^2\psi = tr\big([\mathbf{A}(Du)]D^2\psi\big) = A_{\alpha\beta}(Du)\psi_{,\alpha\beta},$$
$$\text{where } A_{\alpha\beta}(Du) := \mathbf{e}_\alpha \cdot \big(\mathbf{A}(Du)\mathbf{e}_\alpha\big), \tag{3.29}$$

for all smooth scalar fields $\psi$ on $S^2$, and where $\psi_{,\alpha\beta}$ is as defined immediately after (3.8). It is advantageous at this point to expand the second-order terms in (3.28) in order to identify the principal parts:

$$\mathbf{A}(Du) \cdot \left[ D^2\big(B(\phi)H\big) - 2HD^2\big(E(\phi)\big)\right]$$
$$+ \exp(-u)(1 + |Du|^2)^{1/2}\big(\mathbf{A}(Du)[D^2u - Du \otimes Du - \mathbf{1}_x]\mathbf{A}(Du)\big) \cdot D^2\big(E(\phi)\big) \tag{3.30}$$
$$= \mathbf{M}(D^2u, Du, u, \phi) \cdot D^2\phi + B(\phi)\mathbf{A}(Du) \cdot D^2H + ...,$$

where

$$\mathbf{M}(D^2u, Du, u, \phi) := E'(\phi)\exp(-u)(1 + |Du|^2)^{1/2}\big(\mathbf{A}(Du)[D^2u - Du \otimes Du - \mathbf{1}_x]\mathbf{A}(Du)\big)$$
$$+ H[B'(\phi) - 2E'(\phi)]\mathbf{A}(Du).$$

We now substitute (3.9) and (3.10) into (3.27) and (3.28) to obtain the final form of the equilibrium equations that we use for analysis. While the latter leads to a lengthy expression, it is sufficient for our purposes here to identify the principal parts, revealing a quasi-linear system. In particular, the second term in (3.29) becomes

$$B(\phi)\mathbf{A}(Du) \cdot D^2H = \frac{1}{2}\exp(-u)(1 + |Du|^2)^{-3/2} B(\phi)\mathbf{A}(Du) \cdot \big(D^4u[\mathbf{A}(Du)]\big) + ... , \tag{3.31}$$

where we interpret $D^4u(\mathbf{x}) \in L\big(L(T_x, T_x), L(T_x, T_x)\big)$ in (3.31). In components we have

$$\mathbf{A}(Du) \cdot \left( D^4 u[\mathbf{A}(Du)] \right) = A_{\alpha\beta}(Du) A_{\gamma\delta}(Du) u_{,\alpha\beta\gamma\delta},$$
$$\text{where } u_{,\alpha\beta\gamma\delta} := \mathbf{e}_{\alpha} \cdot \left( D^4 \mathbf{u}[\mathbf{e}_{\delta}, \mathbf{e}_{\gamma}, \mathbf{e}_{\beta}] \right), \tag{3.32}$$

the latter notation in accordance with the tri-linear definition of $D^4 u(\mathbf{x})$, cf. (1.9). In view of (3.31), if we multiply (3.28) by the factor $2 \exp(u)(1 + |Du|^2)^{3/2}$, then the equilibrium equations (3.27) and (3.28) take the quasi-linear form

$$-\varepsilon \mathbf{A}(Du) \cdot D^2 \phi + f_1(D^2 u, Du, u, D\phi, \phi, \mu) = 0 \text{ on } S^2, \tag{3.33}$$

and

$$\mathbf{N}(D^2 u, Du, u, \phi) \cdot D^2 \phi + B(\phi) \mathbf{A}(Du) \cdot \left( D^4 u[\mathbf{A}(Du)] \right)$$
$$+ f_2(D^3 u, D^2 u, Du, u, D\phi, \phi, \gamma, \mu) = 0 \text{ on } S^2, \tag{3.34}$$

respectively, where $f_1(\cdot)$ and $f_2(\cdot)$ are $C^2$ functions of their arguments, limited by the smoothness of $B''(\cdot)$ and $E''(\cdot)$, and

$$\mathbf{N}(D^2 u, Du, u, \phi) := 2 \exp(u)(1 + |Du|^2)^{3/2} \mathbf{M}(D^2 u, Du, u, \phi), \tag{3.35}$$

cf. (3.29).

Our problem possesses an obvious "trivial", spherically symmetric solution

$$\phi \equiv \lambda \text{ and } u \equiv 0. \tag{3.36}$$

Using (3.9), (3.10), (3.11), (3.14) and (3.16), we find the remainder of the equilibrium conditions:

$$H \equiv -1, \ K \equiv 1, \ \mu = \Psi'(\lambda),$$
$$\gamma = W(\lambda) - \lambda \Psi'(\lambda) - p / 2, \tag{3.37}$$
$$\Psi(\lambda) := W(\lambda) + B(\lambda) + E(\lambda).$$

Accordingly, we introduce new variables, characterizing non-spherical states,

$$\varphi := \phi - \lambda,$$
$$\varsigma := \gamma - [W(\lambda) - \lambda \Psi'(\lambda) - p / 2], \tag{3.38}$$
$$\xi := \mu - \Psi'(\lambda).$$

which we employ throughout the remainder after substitution into (3.33) and (3.34). Also, using (3.6)$_3$ and (3.38), the constraint equations (2.2), expressed in terms of integrals over $S^2$, now have the form

$$\int_{S^2} \exp(2u)(1 + |Du|^2)^{1/2} ds - 4\pi = 0,$$
$$\int_{S^2} \varphi \exp(2u)(1 + |Du|^2)^{1/2} ds = 0. \tag{3.39}$$

Equations (3.33) and (3.34) comprise a quasilinear system for the pair $(\varphi, u)$ on $S^2$, containing the unknowns (multipliers) $\varsigma, \xi \in \mathbb{R}$, cf. (3.38). The constraint equations (3.39) complete our equilibrium field equations. We choose $\lambda \in \mathbb{R}$ as our bifurcation parameter, which now appears explicitly in (3.33), (3.34), courtesy of (3.38). The non-negative internal pressure $p$ is considered fixed but otherwise arbitrary. Our field equations are parametrized by the constitutive functions $B(\cdot), E(\cdot), W(\cdot)$ and the material constant $\varepsilon$, cf. (2.1)-(2.3) and Figure 2.1.

Before closing this section, we make a simple but important observation. Let $C^{k,\alpha}(S^2)$ denote the space of $k$-times Hölder continuously differentiable functions with exponent $\alpha \in (0,1)$.

**Proposition 3.1** The quasilinear system (3.33) and (3.34) is *elliptic* in the following sense: For fixed $u \in C^{2,\alpha}(S^2)$ and $\varphi \in C^{\alpha}(S^2)$, consider a linear system having principal parts

$$
\begin{aligned}
&-\varepsilon \mathbf{A}(Du) \cdot D^2 h_1 + \ldots = 0, \\
&\mathbf{N}(D^2 u, Du, u, \lambda + \varphi) \cdot D^2 h_1 + B(\lambda + \varphi)\mathbf{A}(Du) \cdot \left(D^4 h_2 [\mathbf{A}(Du)]\right) + \ldots = 0 \text{ on } S^2,
\end{aligned}
\tag{3.40}
$$

where any lower-order terms are presumed smooth. Then system (3.40) is elliptic in the sense of [11].

**Proof.** With $u \in C^{2,\alpha}(S^2)$ and $\varphi \in C^{\alpha}(S^2)$, the resulting coefficient functions (on $S^2$) in the principal part of the system (3.40) are of class $C^{\alpha}$ on $S^2$. For notational convenience, we drop their explicit dependence upon $\varphi, u$ and its derivatives. Following [11], we choose weights $s_1 = s_2 = 0$, $t_1 = 2$ and $t_2 = 4$. Then, using (3.29) and (3.32), we consider the $2 \times 2$ determinant

$$
\begin{vmatrix}
-\varepsilon A_{\alpha\beta}\xi_\alpha\xi_\beta & 0 \\
N_{\alpha\beta}\xi_\alpha\xi_\beta & BA_{\alpha\beta}A_{\gamma\delta}\xi_\alpha\xi_\beta\xi_\gamma\xi_\delta
\end{vmatrix} = -\varepsilon B(A_{\alpha\beta}\xi_\alpha\xi_\beta)^3 < 0 \text{ on } S^2,
$$

for all nonzero $(\xi_1, \xi_2) \in \mathbb{R}^2$. The strict inequality follows from the positivity of the symmetric transformation $\mathbf{A}(Du) \in L(T_x, T_x)$, cf. (3.7)$_3$ and the discussion after (3.10). In particular, the non-vanishing of the above determinant gives the result. □

## 4. Abstract Formulation

Since we intend to carry out a global bifurcation analysis, we want to show that (3.33), (3.34) and (3.39) are equivalent to finding the zeros of a compact vector field on an appropriate function space. Henceforth we write $X^{k,\alpha}, k = 0,1,2,\ldots (X^{\alpha} := X^{0,\alpha})$ for the Hölder function spaces when equipped with the usual Hölder norms on $S^2$, denoted $\|\cdot\|_{k,\alpha}$.

The left sides of (3.33), (3.34), using (3.38), together with the left sides of the constraints (3.39) define mappings

$$
\begin{aligned}
&F_1 : \mathbb{R} \times X \times \mathbb{R}^2 \to Y, \ F_2 : X \to \mathbb{R}^2, \\
&X := X^{2,\alpha} \times X^{4,\alpha}, \ Y := X^{\alpha} \times X^{\alpha}.
\end{aligned}
\tag{4.1}
$$

Now define

$$
\begin{aligned}
&F := (F_1, F_2) : \mathbb{R} \times X \times \mathbb{R}^2 \to Y \times \mathbb{R}^2, \text{ and} \\
&w := (\varphi, u), \tau := (\varsigma, \xi), v := (w, \tau) \equiv (\varphi, u, \varsigma, \xi).
\end{aligned}
\tag{4.2}
$$

Then our equilibrium field equations are equivalent to the abstract operator equation

$$F(\lambda, v) \equiv F(\lambda, w, \tau) \equiv (F_1(\lambda, w, \tau), F_2(w)) = 0. \tag{4.3}$$

In view of (3.37)-(3.38), observe that

$$F(\lambda, 0) \equiv 0, \tag{4.4}$$

i.e., we have the *trivial line of solutions*, $v = 0$ for all $\lambda \in \mathbb{R}$, representing the spherically symmetrical state. Here $F_1(\cdot)$ corresponds to the quasi-linear partial differential operators given by the left sides of (3.33), (3.34), while $F_2(\cdot)$ represents the left side of the two constraint equations (3.39). We note that both $X \times \mathbb{R}^2$ and $Y \times \mathbb{R}^2$ are Banach spaces with the obvious product topology. We further observe that $W(\cdot)$ of class $C^3$ and $B(\cdot)$, $E(\cdot)$ of class $C^4$, cf. Section 2 and (3.34), insure that $F(\cdot)$ is a $C^1$ mapping.

Next we decompose the mapping $F_1(\cdot)$, capturing the quasi-linear form inherent in (3.33), (3.34):

$$F_1(w, \tau) = T_{\lambda, w}[w] + \Phi_1(w, \tau), \tag{4.5}$$

where

$$
\begin{aligned}
&T_{\lambda, w}[w] := (R_u[\varphi], N_{\lambda, w}[\varphi] + S_{\lambda, w}[u]), \\
&R_u[\varphi] := -\varepsilon \mathbf{A}(Du) \cdot D^2\varphi, \\
&N_{\lambda, w}[\varphi] := \mathbf{N}(D^2u, Du, u, \lambda + \varphi) \cdot D^2\varphi, \\
&S_{\lambda, w}[u] := B(\lambda + \varphi)\mathbf{A}(Du) \cdot \left(D^4 u[\mathbf{A}(Du)]\right) \text{ on } S^2,
\end{aligned}
\tag{4.6}
$$

and $\Phi_1 : \mathbb{R} \times X \times \mathbb{R}^2 \to Y$ is defined by the lower-order derivative terms, viz., $\Phi_1$ is the composition operator

$$
\begin{aligned}
\Phi_1(\lambda, w, \tau) := \big( & f_1(D^2u, Du, u, D\varphi, \lambda + \varphi, \xi + \Psi'(\lambda)), \\
& f_2(D^3u, D^2u, Du, u, D\varphi, \lambda + \varphi, \varsigma + W(\lambda) - \lambda\Psi'(\lambda) - p/2, \xi + \Psi'(\lambda)).
\end{aligned}
\tag{4.7}
$$

Accordingly $\Phi_1(\cdot)$ is compact by embedding. The full mapping $F(\cdot)$ in (4.3) then has the quasi-linear form

$$F(\lambda, w, \tau) = \left(T_{\lambda, w}, 0\right)[w, \tau] + \Phi(\lambda, w, \tau), \tag{4.8}$$

where the "0" in the principal part of the operator is simply the $2 \times 2$ zero matrix, and

$$\Phi(\lambda, w, \tau) := \left(\Phi_1(\lambda, w, \tau), F_2(\lambda, w)\right). \tag{4.9}$$

Observing that the constraints equations (3.39) only involve lower-order terms as well, it follow that $\Phi : \mathbb{R} \times X \times \mathbb{R}^2 \to Y \times \mathbb{R}^2$ is compact.

For any constant $M > 0$, define the closed ball

$$\mathcal{B}_M := \{ w = (\varphi, u) \in X : \|w\|_X = \|\varphi\|_{2,\alpha} + \|u\|_{4,\alpha} \le M \}. \tag{4.10}$$

For fixed $(\lambda, w) \in \mathbb{R} \times X$ and constant $a \in \mathbb{R}$, we call

$$\left(T_{\lambda,w}, 0\right) + a := \left(T_{\lambda,w}, 0\right) + a\mathbf{1} : X \times \mathbb{R}^2 \to Y \times \mathbb{R}^2, \tag{4.11}$$

with "$\mathbf{1}$" denoting the identity, the *shifted* linear operator. We henceforth drop the "$\mathbf{1}$" as indicated in (4.11).

**Theorem 4.1** *For any* $w = (\varphi, u) \in \mathcal{B}_M$, $\lambda \in \mathbb{R}$, *and for some sufficiently large positive constant* $a \in \mathbb{R}$, *the shifted operator* $\left(T_{\lambda,w}, 0\right) + a : X \times \mathbb{R}^2 \to Y \times \mathbb{R}^2$, *is bijective.*

We carry out the proof of Theorem 4.1 via two lemmas that follow. First, we note that the linear maps $R_u : X^{2,\alpha} \to X^\alpha$ and $S_{\lambda,w} : X^{4,\alpha} \to X^\alpha$ define uniformly elliptic operators – second and fourth order, respectively. In particular, from (4.6) and (4.10), there are positive constants $c$ and $K$, depending on $M$, such that

$$c^{-1}(\xi_1^2 + \xi_2^2) \le A_{\gamma\delta}(Du)\xi_\gamma\xi_\delta \le c(\xi_1^2 + \xi_2^2) \text{ on } S^2 \text{ for all } (\xi_1, \xi_2) \in \mathbb{R}^2,$$
$$\left\| B(\lambda + \varphi) \right\|_\alpha, \left\| A_{\gamma\delta}(Du) \right\|_\alpha \le K, \ \gamma, \delta = 1, 2. \tag{4.12}$$

**Lemma 4.2** *For every* $w = (\varphi, u) \in \mathcal{B}_M$, $\lambda \in \mathbb{R}$, *the operator* $T_{\lambda,w} + a : X \to Y$ *is injective for* $a \in \mathbb{R}$ *sufficiently large.*

**Proof.** We claim that there are positive, real constants $\delta_\beta, C_\beta, \Lambda_\beta$, depending on $\alpha, c, K, M$, cf. (4.12), but independent of $\mu, \psi, h, \lambda$ and $w$, such that

$$\left\| \psi \right\|_{2,\alpha} \le C_1 \left| \mu \right|^{\alpha/2} \left\| (R_u + \mu)[\psi] \right\|_\alpha,$$
$$\left\| h \right\|_{4,\alpha} \le C_2 \left| \mu \right|^{\alpha/4} \left\| (S_w + \mu)[h] \right\|_\alpha, \tag{4.13}$$

for all $(\psi, h) \in X$, and for all $\mu \in \mathbb{C}$ satisfying $\left| \arg \mu \right| \le \pi / 2 - \delta_\beta, \left| \mu \right| \ge \Lambda_\beta, \beta = 1, 2$, where here we presume the usual complexification of the spaces $X^{k,\alpha}$. Inequalities (4.13) are standard Hölder-space versions of Agmon's estimate (in the $L^p$ setting) [1]. We refer to [21], [27], [31], [41], e.g., for details.

Here and throughout the rest of this work, unless stated otherwise, we assume $a \ge \max(\Lambda_1, \Lambda_2)$. By virtue of (4.13), it follows that $R_u + a$ and $S_{\lambda,w} + a$ are each injective. Now consider $(T_{\lambda,w} + a)[(\psi, h)] = 0$, which according to (4.6), is given explicitly by

$$(R_u + a)[\psi] = 0,$$
$$N_{\lambda,w}[\psi] + (S_{\lambda,w} + a)[h] = 0. \tag{4.14}$$

Then $(4.14)_1 \Rightarrow \psi = 0$, in which case $(4.14)_2 \Rightarrow h = 0$, i.e., $T_{\lambda,w} + a : X \to Y$ is injective. $\square$

Next, we consider the one-parameter family of operators

$$T_{\lambda,w,t}[(\psi, h)] := \left( R_u[\psi], (1-t)N_{\lambda,w}[\psi] + S_{\lambda,w}[h] \right), 0 \le t \le 1. \tag{4.15}$$

Observe that $T_{\lambda,w,0} \equiv T_{\lambda,w}$, while $T_{\lambda,w,1}$ is "diagonal". If necessary, we can always adjust the constant in $(4.12)_2$ so that

$$\left\| N_{\gamma\delta}(D^2 u, Du, u, \phi) \right\|_\alpha \le K, \ \gamma, \delta = 1, 2. \tag{4.16}$$

for every $w = (\varphi, u) \in \mathcal{B}_M$, $\lambda \in \mathbb{R}$, where $N_{\gamma\delta}(D^2 u, Du, u, \lambda + \varphi) := \mathbf{e}_\alpha \cdot \left( \mathbf{N}(D^2 u, Du, u, \lambda + \varphi) \mathbf{e}_\alpha \right)$. Note that Proposition 3.1 is valid for the operator $T_{\lambda, w, t}$ for all $t \in [0, 1]$. Hence, with $w = (\varphi, u) \in \mathcal{B}_M$, $\lambda \in \mathbb{R}$, $t \in [0, 1]$, Proposition 3.1, (4.12) and (4.15) insure the uniform Schauder estimate, cf. [11],

$$\left\| (\psi, h) \right\|_X \le C \left\{ \left\| \left( T_{\lambda, w, t} + a \right)[\psi, h] \right\|_Y + \left\| (\psi, h) \right\|_Y \right\}, \tag{4.17}$$

for all $(\psi, h) \in X$, where the constant $C = C(\alpha, c, K, M)$ is independent of $\psi, h, a, w$ and $t$.

**Lemma 4.3.** *For every* $w = (\varphi, u) \in \mathcal{B}_M$, $\lambda \in \mathbb{R}$, $t \in [0, 1]$, *each of the linear operators* $T_{\lambda, w, t} : X \to Y$ *and* $T_{\lambda, w, t} + a : X \to Y$ *is Fredholm of index zero, i.e., the dimension of the null space and the co-dimension of the range are finite-dimensional and equal.*

**Proof.** From (4.17) and Peetre's lemma, cf. [33], [42], we know that $T_{\lambda, w, t} + a$ and $T_{\lambda, w, t}$ are each *semi-Fredholm*, i.e., the null space is finite-dimensional and the range is closed. Since $(\lambda, w, t)$ is connected to $(0, 0, 1)$ in $\mathbb{R} \times \mathcal{B}_M \times [0, t]$, it follows that the operators $T_{\lambda, w, t} + a$ and $T_{0, 0, 1} + a$ have the same Fredholm index, by the continuity of the index [26]. Using $(3.7)_3$, (3.9) and (4.6), we find

$$T_{0, 0, 1} + a = (-\varepsilon\Delta + a, B(0)\Delta^2 + a), \tag{4.18}$$

where $\Delta^2 := \Delta \circ \Delta$. The linear operator (4.18) is easily seen to be bijective for any $a > 0$: For a given $(b_1, b_2) \in Y = X^\alpha \times X^\alpha$, consider $(T_{0, 0, 1} + a)[(\psi, h)] = (b_1, b_2)$, given explicitly

$$\begin{aligned} -\varepsilon\Delta\psi + a\psi &= b_1, \\ B(0)\Delta^2 h + ah &= b_2 \text{ on } S^2. \end{aligned} \tag{4.19}$$

each of which are uniquely solvable on $X^{2, \alpha}$ and $X^{4, \alpha}$, respectively, by well-known arguments, cf. [29]. In particular, the operator (4.18) and thus, $T_{\lambda, w, t} + a$ have Fredholm index zero. Clearly $T_{\lambda, w, t} + a$ and $T_{\lambda, w, t}$ have the same Fredholm index, treating $a \in \mathbb{R}$ as the connecting parameter. $\square$

**Proof of Theorem 4.1.** The result is immediate: By virtue of Lemma 4.3, the linear operator $T_{\lambda, w, 0} + a \equiv T_{\lambda, w} + a : X \to Y$ is Fredholm of index zero, and from Lemma 4.2 it is injective. Hence, it has a closed range with co-dimension zero, and we conclude that $T_{\lambda, w} + a : X \to Y$ is bijective. Finally, in view of (4.11), note that

$$(T_{\lambda, w}, 0) + a \equiv (T_{\lambda, w} + a, aI_2) : X \times \mathbb{R}^2 \to Y \times \mathbb{R}^2, \tag{4.20}$$

where $I_2$ denotes the identity on $\mathbb{R}^2$, and bijectivity is clear. $\square$

We now return to the full mapping (4.8). For any $\lambda \in \mathbb{R}$ and $v = (w, \tau) \equiv (\varphi, u, \varsigma, \xi) \in \mathcal{B}_M \times \mathbb{R}^2$, consider the linear problem

$$\left( \left( T_{\lambda,w}, 0 \right) + a \right) z = -\Phi(\lambda, v) + a v. \tag{4.21}$$

In view of Theorem 4.1, equation (4.21) has a unique solution, $z \in X \times \mathbb{R}^2$, denoted

$$z = \mathcal{C}_{\lambda,w} \left[ -\Phi(\lambda, v) + a v \right] \coloneqq \mathcal{K}(\lambda, v), \tag{4.22}$$

where $\mathcal{C}_w$ denotes the solution operator associated with the left side of (4.21), viz.,

$$\left( \left( T_{\lambda,w}, 0 \right) + a \right) z = \mathrm{b} \Leftrightarrow z = \mathcal{C}_{\lambda,w} \mathrm{b} \text{ for all } \mathrm{b} \in Y \times \mathbb{R}^2. \tag{4.23}$$

For $t = 0$ in (4.15), recall that $T_{\lambda,w} + a : X \to Y$ is bijective, and consequently the second term on the right side inequality (4.17) may be dropped. With that in hand, (4.17), (4.20) and (4.21) yield

$$\left\| z \right\|_{X \times \mathbb{R}^2} \leq C \left[ \left\| \Phi(\lambda, v) \right\|_{Y \times \mathbb{R}^2} + a \left\| v \right\|_{Y \times \mathbb{R}^2} \right], \tag{4.24}$$

where the constant $C = C(\alpha, c, K, M)$ is independent of $z$ and $v$.

Finally, we define

$$E \coloneqq X^{1,\alpha} \times X^{3,\alpha}. \tag{4.25}$$

Then (4.22)-(4.24) and the compact embedding, $E \subset X$, imply that each of the operators $\mathcal{C}_{\lambda,w} :$ $(\mathcal{B}_M \cap E) \times \mathbb{R}^2 \to E \times \mathbb{R}^2$ and $\mathcal{K} : \mathbb{R} \times (\mathcal{B}_M \cap E) \times \mathbb{R}^2 \to E \times \mathbb{R}^2$ is compact. Moreover, (4.22), (4.24) and the continuity of $\Phi(\cdot)$, cf. (4.7), (4.9), imply that $\mathcal{K}(\cdot)$ is continuous, while the same argument applied to (4.23) shows that $w \mapsto \mathcal{C}_w$ is continuous as well in the operator norm associated with $L(Y \times \mathbb{R}^2, E \times \mathbb{R}^2)$. We conclude that (4.4) on $\mathbb{R} \times \mathcal{B}_M \times \mathbb{R}^2$ is equivalent to the operator equation

$$v - \mathcal{K}(\lambda, v) = 0, \tag{4.26}$$

where $\mathcal{K} : \mathbb{R} \times (\mathcal{B}_M \cap E) \times \mathbb{R}^2 \to E \times \mathbb{R}^2$ is continuous and compact. By virtue of (4.4), we note that

$$\mathcal{K}(\lambda, 0) \equiv 0. \tag{4.27}$$

## 5. Linearized Problem

In order to obtain the linearization of the operator (4.2) about the trivial solution, we return to (3.36)-(3.38) and substitute

$$\varphi \equiv \alpha \vartheta, \ u \equiv \alpha \upsilon, \ H = -1 + \alpha h, \ K = 1 + \alpha \kappa,$$
$$\mu = \Psi'(\lambda) + \alpha \xi, \ \gamma = W(\lambda) - \lambda \Psi'(\lambda) - p / 2 + \alpha \varsigma, \tag{5.1}$$

into (3.7)$_3$, (3.9), (3.10), (3.27), (3.28) and (3.39). In each of these we formally compute the directional derivative, $d(\cdot) / d\alpha |_{\alpha=0}$, to deduce the linearized equations:

$$\kappa = -2h,$$
$$\Delta \upsilon + 2\upsilon = 2h,$$
$$-\varepsilon \Delta \vartheta + \Psi''(\lambda)\vartheta = \xi, \tag{5.2}$$
$$B(\lambda)\Delta h - ph + [E'(\lambda) - B'(\lambda)]\Delta \vartheta - 2[B'(\lambda) + E'(\lambda)]\vartheta = 2(\varsigma + \lambda \xi) \text{ on } S^2,$$

$$\int_{S^2} \upsilon ds = 0,$$
$$\int_{S^2} \vartheta ds = 0, \tag{5.3}$$

We first integrate (5.2)$_2$ over $S^2$ and use (1.8), the surface divergence theorem (1.5) and (5.3)$_1$ to deduce

$$\int_{S^2} h ds = 0. \tag{5.4}$$

Applying the same steps to (5.2)$_{3,4}$ while employing (5.3)$_2$ and (5.4), we conclude that

$$\varsigma = \xi = 0. \tag{5.5}$$

In order to proceed, recall that the eigenvalues of the Laplace-Beltrami operator on $S^2$ are given by $-\ell(\ell+1),\ \ell = 0,1,2,...,$ with accompanying eigenfunctions (surface harmonics)

$$\rho_{\ell,m}(\mathbf{x}),\ m = -\ell,...,\ell,$$

where $\rho_{\ell,0}(\mathbf{x}) := P_\ell(\cos\theta),$

$$\rho_{\ell,m}(\mathbf{x}) := P_{\ell,m}(\cos\theta)\cos m\psi,\ m = 1,...,\ell,$$
$$\rho_{\ell,m}(\mathbf{x}) := -P_{\ell,-m}(\cos\theta)\sin m\psi,\ m = -\ell,...,-1,\ l = 1,2,... \tag{5.6}$$

Here we express

$$\mathbf{x} = (x_1, x_2, x_3) = (\sin\theta\cos\psi, \sin\theta\sin\psi, \cos\theta) \in S^2, \tag{5.7}$$

where $\theta \in [0, \pi],\ \psi \in [0, 2\pi)$, denote the usual spherical coordinates, and $P_\ell(\cdot)$ and $P_{\ell,m}(\cdot)$ denote the Legendre and associated Legendre polynomials, respectively.

We consider (5.2)$_3$ and (5.3)$_2$, using (5.5), which are decoupled from the rest of (5.2), (5.3), viz.,

$$-\varepsilon \Delta \vartheta + \Psi''(\lambda)\vartheta = 0 \text{ on } S^2,$$
$$\int_{S^2} \vartheta ds = 0. \tag{5.8}$$

Clearly (5.8) admits nontrivial solutions if and only if $\lambda$ and $\varepsilon$ satisfy the *characteristic equation*

$$-\Psi''(\lambda) / \varepsilon = \ell(\ell+1), \text{ for } \ell \in \mathbb{N}. \tag{5.9}$$

For a given $\ell \in \mathbb{N}$ and $\varepsilon$ sufficiently small, we assume that (5.9) has at least one root, denoted $\lambda_\ell$. For example, in the special case that the bending modulus function $B$ and $E$ are constant, we have $\Psi(\lambda) \equiv W(\lambda)$, cf. (5.2)$_3$. From the graph of $W''(\cdot)$ shown in Figure 2.1, we see that $\lambda \in (m_1, m_2)$ is necessary in order for (5.9) to have roots in this case. Depending on the size of $\varepsilon$, there can be no roots, one root, or precisely two roots, the latter two possibilities corresponding to $\varepsilon$ sufficiently small. Returning to (5.9), it is reasonable to assume that $\Psi(\cdot)$ has properties similar to $W(\cdot)$. In particular, we assume that $\Psi''(\cdot) < 0$ on some open (spinodal) set, with $\varepsilon$ sufficiently small. Given a root of the characteristic equation, we then have the nontrivial solution

$$\vartheta = \vartheta_\ell := \sum_{m=-\ell}^{\ell} c_m^\ell \rho_{\ell,m}, \tag{5.10}$$

for arbitrary constants $c_m^l$.

Next we turn to (5.2)$_2$ and (5.3)$_1$. We first note that the operator on the left side of (5.2)$_2$ has the nontrivial homogeneous solution, corresponding to $\ell = 1$:

$$\upsilon = \upsilon_1 := \sum_{m=-1}^{1} c_m^1 \rho_{1,m}, \tag{5.11}$$

Since the operator on the left side is formally self-adjoint, it follows that (5.2)$_2$ has a particular solutions if and only if the right side is $L^2$- orthogonal to all such $\upsilon_1$ given by (5.11). Accordingly, we consider (5.2)$_4$ and (5.4) at $\lambda = \lambda_\ell$ and $\vartheta = \vartheta_\ell$, for $\ell \geq 2$, cf. (5.10):

$$B(\lambda_\ell)\Delta h - ph = \left( E'(\lambda_\ell)[2 + \ell(\ell+1)] + B'(\lambda_\ell)[2 - \ell(\ell+1)] \right) \varphi_\ell. \tag{5.12}$$

Recall that $p \geq 0$ is imposed, cf. Section 1, in which case the left side of (5.12) admits no homogeneous solutions, cf. (2.3). Hence, for $\ell \geq 2$, (5.12) has the solution

$$
\begin{aligned}
h &= h_\ell := \sigma_\ell \vartheta_\ell = \sigma_\ell \sum_{m=-\ell}^{\ell} c_m^\ell \rho_{\ell,m}, \\
\sigma_\ell &:= -\frac{\left( E'(\lambda_\ell)[2 + \ell(\ell+1)] + B'(\lambda_\ell)[2 - \ell(\ell+1)] \right)}{[\ell(\ell+1)B(\lambda_\ell) + p]}.
\end{aligned}
\tag{5.13}
$$

Returning now to (5.2)$_2$ with $h = h_\ell, \ell \geq 2$, we find the particular solution

$$
\begin{aligned}
\upsilon &= \upsilon_\ell := \tau_\ell \vartheta_\ell = \tau_\ell \sum_{m=-\ell}^{\ell} c_m^\ell \rho_{\ell,m}, \\
\tau_\ell &:= -2 \frac{\left( E'(\lambda_\ell)[2 + \ell(\ell+1)] + B'(\lambda_\ell)[2 - \ell(\ell+1)] \right)}{[\ell(\ell+1)B(\lambda_\ell) + p][2 - \ell(\ell+1)]}.
\end{aligned}
\tag{5.14}
$$

In order to rigorously summarize our efforts here, we define the finite-dimensional subspaces

$$
\begin{aligned}
\mathcal{N}_o &:= span\{(0, \rho_{1,j}, 0, 0), \ j = -1, 0, 1\} \subset X \times \mathbb{R}^2, \\
\mathcal{N}_\ell &:= span\{z_{\ell,m} := (\rho_{\ell,m}, \tau_\ell \rho_{\ell,m}, 0, 0), \ m = -\ell, ..., \ell\} \subset X \times \mathbb{R}^2, \ell \in \mathbb{N},
\end{aligned}
\tag{5.15}
$$

where $\tau_\ell$ is given by (5.14)$_2$ for all $\ell \geq 2$, and $\tau_1 = 0$. In view of (4.2)-(4.4), we have **Proposition 5.1** *The Fréchet derivative of $v \mapsto F(\lambda, v)$ evaluated along the trivial solution, denoted*
$L(\lambda) := D_v F(\lambda, 0) : X \times \mathbb{R}^2 \to Y \times \mathbb{R}^2$, *is given by*

$$
\begin{aligned}
L(\lambda)[z] = \big( &-\varepsilon\Delta\vartheta + \Psi''(\lambda)\vartheta - \xi, 2[E'(\lambda) - B'(\lambda)]\Delta\vartheta + B(\lambda)\Delta^2 \upsilon \\
&+ [2B(\lambda) - p]\Delta\upsilon - 2p\upsilon - 4[B'(\lambda) + E'(\lambda)]\vartheta - 4(\varsigma + \lambda\xi), \\
&\qquad\qquad\qquad\qquad\qquad\qquad\qquad \int_{S^2}\upsilon ds, \int_{S^2}\vartheta ds \big),
\end{aligned}
\tag{5.16}
$$

*for all $z = (\vartheta, \upsilon, \varsigma, \xi) \in X \times \mathbb{R}^2$. Let $\lambda = \lambda_\ell$ denote a of the root of the characteristic equation (5.9) for $\ell \in \mathbb{N}$. The null space of $L(\lambda)$ is characterized by*

$$
\begin{aligned}
\mathcal{N}\big(L(\lambda)\big) &= \mathcal{N}_o, \text{ for all } \lambda \neq \lambda_\ell; \\
\mathcal{N}\big(L(\lambda_\ell)\big) &= \mathcal{N}_\ell \oplus \mathcal{N}_o,
\end{aligned}
\tag{5.17}
$$

*cf. (5.15). Moreover, $L(\lambda) : X \times \mathbb{R}^2 \to Y \times \mathbb{R}^2$ is a Fredholm operator of index zero.*

**Proof.** As previously mentioned, differentiability here is a direct consequence of the presumed smoothness of the functions $W(\cdot), B(\cdot), E(\cdot)$, cf. (3.34). The second component on the right side of (5.16) results from the substitution of $(5.2)_2$ into $(5.2)_4$, keeping in mind (3.31) and the elimination of the factor "1/2" leading to (3.34). The first two claims summarize the calculations (5.1)-(5.14). Finally, note that the principal part of the operator $L(\lambda)$ is simply $\big(T_{\lambda,0}, 0\big)$, i.e., $L(\lambda) = \big(T_{\lambda,0}, 0\big) + ...$, where

$$
\big(T_{\lambda,0}, 0\big)[z] \equiv \big(-\varepsilon\Delta\vartheta, 2[E'(\lambda) - B'(\lambda)]\Delta\vartheta + B(\lambda)\Delta^2 \upsilon, 0, 0\big),
\tag{5.18}
$$

for all $z = (\vartheta, \upsilon, \xi, \varsigma) \in X \times \mathbb{R}^2$. Clearly the two-dimensional extension of $T_{\lambda,0}$ on $X$ to $\big(T_{\lambda,0}, 0\big)$ on $X \times \mathbb{R}^2$ does not change the index. Recalling $T_{\lambda,0} \equiv T_{\lambda,0,0}$, cf. (4.15), the last claim follows from Lemma 4.3 and the stability of the Fredholm index, i.e., $L(\lambda)$ and $\big(T_{\lambda,0}, 0\big)$ have the same zero Fredholm index. $\square$

In order to obtain global bifurcation results, we also need to consider some properties of the continuous, compact mapping in (4.26), cf. (4.21)-(4.25):

**Proposition 5.2.** *The mapping $v \mapsto \mathcal{K}(\lambda, v), \mathcal{K}(\lambda, \cdot) : E \times \mathbb{R}^2 \to E \times \mathbb{R}^2$, is differentiable along the trivial line, $v \equiv 0$, with Fréchet derivative, denoted $\mathcal{A}(\lambda) := D_v \mathcal{K}(\lambda, 0) \in L(E \times \mathbb{R}^2, E \times \mathbb{R}^2)$, given by*

$$
\mathcal{A}(\lambda)[z] = -\mathcal{C}_{\lambda,0}\big([D_v\Phi(\lambda, 0) - a]z\big),
\tag{5.19}
$$

*for all $z = (\vartheta, \upsilon, \varsigma, \xi) \in E$, where $D_v\Phi(\lambda, 0)$ denotes the Fréchet derivative of $v \mapsto \Phi(\lambda, v)$ at $v \equiv 0$, with $\Phi(\lambda, \cdot) : E \times \mathbb{R}^2 \to Y \times \mathbb{R}^2$. Moreover, $\lambda \mapsto \mathcal{A}(\lambda)$ is continuous, and $\mathcal{A}(\lambda)$ is compact.*

**Proof.** Since $\Phi(\lambda, \cdot)$ involves only lower-order derivatives, cf. (4.7)-(4.9), it follows that $\Phi(\cdot)$ is $C^1$ on $E \times \mathbb{R}^2$ as well as on $X \times \mathbb{R}^2$. Moreover, from (4.4) and (4.8), we see that $\Phi(\lambda, 0) \equiv 0$. Denoting $q := (\vartheta, \upsilon)$, while bearing in mind (4.27), consider

$$\left\| \mathcal{K}(\lambda, z) + \mathcal{C}_{\lambda, 0}\left( [D_v \Phi(\lambda, 0) - a] z \right) \right\|$$

$$= \left\| \mathcal{C}_{\lambda, 0}\left( [D_v \Phi(\lambda, 0) - a] z \right) - \mathcal{C}_{\lambda, q}\left( \Phi(\lambda, z) - az \right) \right\|$$

$$\leq \left\| \mathcal{C}_{\lambda, 0} - \mathcal{C}_{\lambda, q} \right\|_L \left\| \Phi(\lambda, z) - az \right\| + \left\| \mathcal{C}_{\lambda, 0} \right\|_L \left\| \Phi(\lambda, z) - D_v \Phi(\lambda, 0) z \right\| \qquad (5.20)$$

$$\leq \left\| \mathcal{C}_{\lambda, 0} - \mathcal{C}_{\lambda, q} \right\|_L \left\| D_v \Phi(\lambda, 0) z - az \right\|$$

$$\qquad + \left( \left\| \mathcal{C}_{\lambda, 0} \right\|_L + \left\| \mathcal{C}_{\lambda, 0} - \mathcal{C}_{\lambda, q} \right\|_L \right) \left\| \Phi(\lambda, z) - D_v \Phi(\lambda, 0) z \right\|,$$

where $\left\| \cdot \right\|_L$ denotes the operator norm on the space $L(Y \times \mathbb{R}^2, E \times \mathbb{R}^2)$. The differentiability of $\Phi(\cdot)$ implies that, in particular, for each $\lambda \in \mathbb{R}$ and a given $\varepsilon > 0$, there is a $\delta_\lambda > 0$ such that $\left\| z \right\| \leq \delta_\lambda \Longrightarrow$

$$\left\| \Phi(\lambda, z) - D_v \Phi(\lambda, 0) z \right\| \leq \varepsilon \left\| z \right\|, \qquad (5.21)$$

while the continuity of $(\lambda, q) \mapsto \mathcal{C}_{\lambda, q}$ implies, say,

$$\left\| \mathcal{C}_{\lambda, 0} - \mathcal{C}_{\lambda, q} \right\|_L \leq \varepsilon. \qquad (5.22)$$

We also deduce

$$\left\| D_v \Phi(\lambda, 0) z - az \right\| \leq \left\| D_v \Phi(\lambda, 0) - a \right\|_{\tilde{L}} \left\| z \right\| \coloneqq M(\lambda) \left\| z \right\|, \qquad (5.23)$$

where $M(\lambda) < \infty$, and $\left\| \cdot \right\|_{\tilde{L}}$ denotes the operator norm on $L(E \times \mathbb{R}^2, Y \times \mathbb{R}^2)$. With (5.21)-(5.23) in hand, for a given $\varepsilon > 0$, the right side of inequality (5.20), is bounded above by $\varepsilon C(\lambda) \left\| z \right\|$ for all $\left\| z \right\| \leq \delta_\lambda$, where $C(\lambda) < \infty$. Lastly, (5.19) and the continuity of the maps $\lambda \mapsto \mathcal{C}_{\lambda, 0}, D_v \Phi(\lambda, 0)$ show that $\mathcal{A}(\cdot)$ is continuous, while the derivative of a compact map is compact. □

We close this section with

**Proposition 5.3.** *Let $\mathcal{I}$ denote the identity on $E \times \mathbb{R}^2$. Then the null space of $\mathcal{I} - \mathcal{A}(\lambda)$ is the same as that of $L(\lambda)$ as described in* (5.17), *i.e.,* $\mathcal{N}(\mathcal{I} - \mathcal{A}(\lambda)) = \mathcal{N}(L(\lambda))$.

**Proof.** This follows easily from (4.21), (4.23), (5.16), (5.18) and (5.19):

$$L(\lambda)[z] = (T_{\lambda, 0}, 0)[z] + D_v \Phi(\lambda, 0)[z]$$

$$= \left( (T_{\lambda, 0}, 0) + a \right)[z] + \left( D_v \Phi(\lambda, 0) - a \right)[z] = 0 \qquad (5.24)$$

$$\Leftrightarrow z + \mathcal{C}_{\lambda, 0}\left( D_v \Phi(\lambda, 0) - a \right)[z] = z - \mathcal{A}(\lambda)[z] = 0. \ \square$$

## 6. Global Symmetry-Breaking Bifurcation

We begin this section by identifying the equivariant symmetries of the problem embodied in (4.3). Let $O(3)$ denote the orthogonal group acting on $\mathbb{R}^3$. We define the *natural* action of $\mathbf{G} \in O(3)$ on any $v = (\varphi, u, \varsigma, \xi) \in Y \times \mathbb{R}^2$, via

$$\Gamma v \coloneqq (\varphi(\mathbf{G}^T \mathbf{x}), u(\mathbf{G}^T \mathbf{x}), \varsigma, \xi) \text{ for all } \mathbf{G} \in O(3). \qquad (6.1)$$

Referring to the mapping $F(\cdot)$, cf. (4.1)-(4.3), we claim

**Proposition 6.1.** *The nonlinear operator $F(\lambda,\cdot)$ is equivariant under the natural action* (6.1):

$$F(\lambda,\Gamma v) = \Gamma F(\lambda, v) \text{ for all } \mathbf{G} \in O(3). \tag{6.2}$$

**Proof.** The proof, which follows along the lines of [18], is straightforward but tedious. Here we indicate some of the necessary steps. First, the transformations

$$D\varphi(\mathbf{x}), Du(\mathbf{x}) \to \mathbf{G}D\varphi(\mathbf{G}^T\mathbf{x}), \mathbf{G}Du(\mathbf{G}^T\mathbf{x}),$$

*and* $\tag{6.3}$

$$D^2\varphi(\mathbf{x}), D^2u(\mathbf{x}) \to \mathbf{G}[D^2\phi(\mathbf{G}^T\mathbf{x})]\mathbf{G}^T, \mathbf{G}[D^2u(\mathbf{G}^T\mathbf{x})]\mathbf{G}^T, \text{ for all } \mathbf{G} \in O(3),$$

easily follow from $v \to \Gamma v$ via (6.1). From (6.2) and a bit of work, we also find

$$
\begin{aligned}
&H(\mathbf{x}) \to H(\mathbf{G}^T\mathbf{x}),\\
&K(\mathbf{x}) \to K(\mathbf{G}^T\mathbf{x}),\\
&\mathbf{A}(Du(\mathbf{x})) \to \mathbf{G}A(Du(\mathbf{G}^T\mathbf{x}))\mathbf{G}^T,\\
&g(D^2u(\mathbf{x}), Du(\mathbf{x})) \to g(D^2u(\mathbf{G}^T\mathbf{x}), Du(\mathbf{G}^T\mathbf{x})),\\
&\mathbf{r}(D^2u(\mathbf{x}), Du(\mathbf{x}), u(\mathbf{x})) \to \mathbf{G}\mathbf{r}(D^2u(\mathbf{G}^T\mathbf{x}), Du(\mathbf{G}^T\mathbf{x}), u(\mathbf{G}^T\mathbf{x})), \text{ for all } \mathbf{G} \in O(3),
\end{aligned}
\tag{6.4}
$$

where the quantities involved in (6.4) are defined in (3.9), (3.10), $(3.7)_3$, (3.20) and (3.26), respectively.

We give a sampling of the calculations involved leading to (6.4): Recall from (3.2) that $\mathbf{1}_x \equiv \mathbf{I} - \mathbf{x} \otimes \mathbf{x}$, is the identity on the tangent space $T_x$, where $\mathbf{x}$ denotes the unit normal, interpreted here as the unit-radial translation of $\mathbf{x} \in S^2$. Then from $(3.7)_3$ and (6.3) we then have

$$
\begin{aligned}
\mathbf{A}(Du(\mathbf{x})) &\to (1+|\mathbf{G}Du(\mathbf{G}^T\mathbf{x})|^2)\{\mathbf{1}_x - (\mathbf{G}Du(\mathbf{G}^T\mathbf{x})) \otimes (\mathbf{G}Du(\mathbf{G}^T\mathbf{x}))\}\\
&= (1+|Du(\mathbf{G}^T\mathbf{x})|^2)\{\mathbf{G}[\mathbf{I} - \mathbf{G}^T\mathbf{x} \otimes \mathbf{G}^T\mathbf{x}]\mathbf{G}^T - \mathbf{G}[Du(\mathbf{G}^T\mathbf{x}) \otimes Du(\mathbf{G}^T\mathbf{x})]\mathbf{G}^T\}\\
&= (1+|Du(\mathbf{G}^T\mathbf{x})|^2)\mathbf{G}\{\mathbf{1}_{\mathbf{G}^T x} - Du(\mathbf{G}^T\mathbf{x}) \otimes Du(\mathbf{G}^T\mathbf{x})\}\mathbf{G}^T,\\
&\equiv \mathbf{G}A(Du(\mathbf{G}^T\mathbf{x}))\mathbf{G}^T,
\end{aligned}
\tag{6.5}
$$

where we have used the identity $\mathbf{a} \otimes \mathbf{G}^T\mathbf{b} = \mathbf{a} \otimes \mathbf{b}\mathbf{G}$. From (3.9), (6.3) and (6.5) we also find

$$
\begin{aligned}
H(\mathbf{x}) &\to \exp(-u(\mathbf{G}^T\mathbf{x}))(1+|\mathbf{G}Du(\mathbf{G}^T\mathbf{x})|^2)^{-3/2}\{(\mathbf{G}[\mathbf{A}(\mathbf{G}^T\mathbf{x})]\mathbf{G}^T) \cdot (\mathbf{G}[D^2u(\mathbf{G}^T\mathbf{x})]\mathbf{G}^T)\\
&\qquad\qquad\qquad - 2(1+|\mathbf{G}Du(\mathbf{G}^T\mathbf{x})|^2)\}/2\\
&= \exp(-u(\mathbf{G}^T\mathbf{x}))(1+|Du(\mathbf{G}^T\mathbf{x})|^2)^{-3/2}\{(\mathbf{A}(\mathbf{G}^T\mathbf{x})) \cdot (D^2u(\mathbf{G}^T\mathbf{x}))\\
&\qquad\qquad\qquad - 2(1+|Du(\mathbf{G}^T\mathbf{x})|^2)\}/2\\
&\equiv H(\mathbf{G}^T\mathbf{x}).
\end{aligned}
\tag{6.6}
$$

Similar calculations for (3.12), (3.21) and (3.22) (the building blocks for $g(\cdot)$ and $\mathbf{r}(\cdot)$) show that

$$
\begin{aligned}
\mathbf{F}^{-1}(\mathbf{x}) &\to \mathbf{G}[\mathbf{F}^{-1}(\mathbf{G}^T\mathbf{x})]\mathbf{G}^T, \\
\mathbf{P}_y(\mathbf{x}) &\to \mathbf{G}[\mathbf{P}_y(\mathbf{G}^T\mathbf{x})]\mathbf{G}^T, \\
\mathbf{L}(\mathbf{x}) &\to \mathbf{G}[\mathbf{L}(\mathbf{G}^T\mathbf{x})]\mathbf{G}^T, \\
\nabla\mathbf{F}^{-T}(\mathbf{x})\nabla\varphi(\mathbf{x}) &\to \mathbf{G}[\nabla\mathbf{F}^{-T}(\mathbf{G}^T\mathbf{x})\nabla\varphi(\mathbf{G}^T\mathbf{x})]\mathbf{G}^T.
\end{aligned}
\tag{6.7}
$$

The definition (3.20) and the transformations $(6.7)_{1,2,4}$ give $(6.4)_4$, and likewise (3.26) along with all in (6.7) yield $(6.4)_5$.

With (6.3), (6.4) and $(6.7)_3$ in hand, it is straightforward to verify that both (3.27) (hence (3.33)) and (3.28) transform like the first two components of (6.1), while the constraint equations (3.39) are invariant - like the last two components of (6.1). Finally, the equivariance of (3.34) in terms of the natural action on $u$, follows by composing (3.28) with (3.9) and then using (6.6) and the equivariance of (3.28). Note that the scalar factor in (3.31), involved in the final form of (3.34), is invariant and thus has no effect on equivariance of the latter. □

We suppose now that $\lambda_\ell$ satisfies the characteristic equation (5.9) for a given $\ell \in \mathbb{N}$ with $\varepsilon$ sufficiently small, i.e., $(\lambda_\ell, 0) \in \mathbb{R} \times X \times \mathbb{R}^2$ is a potential bifurcation point. Aside from the 3-dimensional space $\mathcal{N}_o$, the high dimensionality of the null space, cf. (5.15) and (5.17), is a direct consequence of the $O(3)$ symmetry; from (6.2) the linear operator (5.16) also commutes with the action (6.1). As is well known, the equivariance properties (6.2) can be used to simplify an otherwise difficult bifurcation analysis, cf. [5], [14], [36], [40]. In order to obtain global results in a straightforward manner, we proceed as in [18]: For any subgroup $\mathcal{G} \subset O(3)$, we define the fixed-point spaces

$$
\begin{aligned}
X_\mathcal{G} &:= \{w \in X : \Gamma w = w \text{ for all } \mathbf{G} \in \mathcal{G}\}, \\
Y_\mathcal{G} &:= \{y \in Y : \Gamma y = y \text{ for all } \mathbf{G} \in \mathcal{G}\},
\end{aligned}
\tag{6.8}
$$

which are each closed subspaces. By virtue of (6.2) we then have

$$
F : \mathbb{R} \times X_\mathcal{G} \times \mathbb{R}^2 \to Y_\mathcal{G} \times \mathbb{R}^2,
\tag{6.9}
$$

i.e., the nonlinear mapping $F(\lambda, \cdot)$ has linear invariant subspaces.

The basic idea then is to strategically choose subgroups $\mathcal{G} \subset O(3)$ so that a bifurcation analysis of (4.3) or (4.26) can be carried out easily on the reduced space via (6.9). This is particularly the case for any subgroup $\mathcal{G}$ leading to a one-dimensional null space, viz.,

$$
\dim \mathcal{N}\left(L(\lambda_\ell)\big|_{X_\mathcal{G} \times \mathbb{R}^2}\right) = 1.
\tag{6.10}
$$

If (6.10) holds along with the crossing condition,

$$
L'(\lambda_\ell)[z] \notin \mathcal{R}\left(L(\lambda_\ell)\big|_{X_\mathcal{G} \times \mathbb{R}^2}\right), \ z \in \mathcal{N}\left(L(\lambda_\ell)\big|_{X_\mathcal{G} \times \mathbb{R}^2}\right),
\tag{6.11}
$$

then (4.3) has a local branch of bifurcating solutions in $\mathbb{R} \times X_\mathcal{G} \times \mathbb{R}^2$. Assuming that (6.10) is true, then the condition

$$\Psi'''(\lambda_\ell) \neq 0, \tag{6.12}$$

insures that (6.11) is fulfilled, as we explain below. We now summarize:

**Proposition 6.2** *Assume that $\lambda_\ell$ is a root of the characteristic equation (5.9) for a given $\ell \in \mathbb{N}$, such that*

$$\mathcal{N}\left(L(\lambda_\ell)\big|_{X_\mathcal{G} \times \mathbb{R}^2}\right) = span\{\hat{z}_\ell\}, \tag{6.13}$$

*for $\hat{z}_\ell \in \mathcal{N}_\ell$, with $X_\mathcal{G}$ is defined in (6.8) for some subgroup $\mathcal{G} \subset O(3)$. If (6.12) holds, then there is a unique, $C^1$ curve of local, nontrivial solution pairs for (4.3):*

$$\begin{aligned}
&\{(\lambda, v) = (\hat{\lambda}(t), \hat{v}(t)) : |t| < \delta\} \subset \mathbb{R} \times X_\mathcal{G} \times \mathbb{R}^2, \\
&F(\hat{\lambda}(t), \hat{v}(t)) \equiv 0, \\
&\hat{\lambda}(t) = \lambda_\ell + O(t), \ \hat{v}(t) = t\hat{z}_\ell + o(t), \ as \ t \to 0,
\end{aligned} \tag{6.14}$$

*In addition, if $\ell \in \mathbb{N}$ is odd, then $\hat{\lambda}(t) = \lambda_\ell + o(t)$, i.e., the solution branch is a "pitchfork".*

**Proof.** The proof, based on a local Liapunov-Schmidt reduction, is standard, e.g., [28]. We merely indicate the role of (6.12). If the latter holds, then from (5.15), (5.16), the first component of the left side of (6.11) equals $\Psi'''(\lambda_\ell)\hat{\rho}_\ell$, where $\hat{z}_\ell := (\hat{\rho}_\ell, \tau_\ell \hat{\rho}_\ell, 0, 0)$. Using the first and fourth components of (5.16), it is clear that a term with first component proportional to an eigenfunction, $\hat{\rho}_\ell$, is not in the range of the operator at $\lambda_\ell$, as required in (6.11). That $\ell \in \mathbb{N}$ odd gives rise to "pitchforks" while $\ell \in \mathbb{N}$ even leads to "transcritical" local bifurcations is well known. □

Recall that any solution of (4.26) is also a solution of (4.3) and vice-versa. Thus, (6.13) yields a local nontrivial solution curve to (4.26) as well. Of course we could also demonstrate that the compact vector field on the left side of (4.26) is equivariant. But it is more expedient to realize that all of Section 4 is valid with $X_\mathcal{G}$ and $Y_\mathcal{G}$ in place of $X$ and $Y$, respectively, again resulting in (4.26), with

$\mathcal{K} : \mathbb{R} \times (\mathcal{B}_M \cap E_\mathcal{G}) \times \mathbb{R}^2 \to E_\mathcal{G} \times \mathbb{R}^2$, where $E_\mathcal{G} := E \cap X_\mathcal{G}$. This is the starting point for extending Proposition 6.1 to global conclusions.

In view of (6.13), Propositions 5.1 and 5.3 imply that both operators, $L(\lambda_\ell)\big|_{X_\mathcal{G} \times \mathbb{R}^2}$ and

$[\mathcal{I} - \mathcal{A}(\lambda_\ell)]\big|_{E_\mathcal{G} \times \mathbb{R}^2}$, have an isolated zero eigenvalue of finite algebraic multiplicity. In fact, the zero eigenvalue here is simple, as readily shown by direct calculation (as in Section 5). A well-known perturbation argument, via the implicit function theorem [28], then insures the existence of a local, differentiable eigenvalue curve, say, $\sigma = \hat{\sigma}(\lambda)$, with $\hat{\sigma}(\lambda_\ell) = 0$, for the eigenvalue problem

$$L(\lambda)[z] = \sigma z, \ z \in X_\mathcal{G} \times \mathbb{R}^2. \tag{6.15}$$

Fortunately, the first and last components of (6.15) decouple from the rest, cf. (5.16), leading to

$$\begin{aligned}
&-\varepsilon \Delta \vartheta + \Psi''(\lambda)\vartheta - \xi = \sigma \vartheta, \\
&\int_{S^2} \vartheta ds = \sigma \xi.
\end{aligned}$$

This yields the explicit eigenvalue-perturbation curve

$$\sigma = \varepsilon\ell(\ell+1) + \Psi''(\lambda) \coloneqq \hat\sigma(\lambda), \tag{6.16}$$

with $\vartheta = \hat\rho_\ell$ and $\xi = 0$, which verifies $\hat\sigma(\lambda_\ell) = 0$ because $\lambda_\ell$ is a root of the characteristic equation (5.9). We remark that the other components comprising the eigenvector, say, $z = z_\lambda$ in (6.15) now follow by back-substitution, the precise form of which is not important for our purposes here. In any case, we now see from (6.16) that

$$\hat\sigma'(\lambda_\ell) = \Psi'''(\lambda_\ell), \tag{6.17}$$

i.e., (6.12) implies that $\hat\sigma(\lambda)$ has a simple zero at $\lambda_\ell$.

We now rewrite (6.15) as in (5.24), at $(\sigma, z) = (\hat\sigma(\lambda), z_\lambda)$, to deduce

$$z_\lambda - \mathcal{A}(\lambda)[z_\lambda] = \hat\sigma(\lambda)\mathcal{C}_{\lambda,0}[z_\lambda]. \tag{6.18}$$

In order to understand the right side of (6.18), we consider the eigenvalue problem

$$\big((T_{\lambda,0}, 0) + a\big)[z] = \omega z, \; z \in X_{\mathcal{G}} \times \mathbb{R}^2. \tag{6.19}$$

From (5.18), the first component of (6.19) decouples from the rest, and for $z = z_\lambda$, we find that $\omega = a + \varepsilon\ell(\ell+1)$, which is independent of $\lambda$. We then use (4.23) to deduce

$$z_\lambda = \big(a + \varepsilon\ell(\ell+1)\big)\mathcal{C}_{\lambda,0}[z_\lambda]. \tag{6.20}$$

Incorporating this into (6.18) leads to the eigenvalue perturbation

$$z_\lambda - \mathcal{A}(\lambda)[z_\lambda] = \frac{\hat\sigma(\lambda)}{a + \varepsilon\ell(\ell+1)}\, z_\lambda. \tag{6.21}$$

We may now state a global bifurcation result [35]:

**Theorem 6.3.** *Assume the hypotheses of Proposition 6.2 concerning* (6.13)*. In addition, suppose that the function*

$$\hat\sigma(\lambda) = \varepsilon\ell(\ell+1) + \Psi''(\lambda) \text{ changes sign at } \lambda = \lambda_\ell. \tag{6.22}$$

*Let $\mathcal{S}$ denote the closure of all nontrivial solutions of* (4.26)*. Then $(\lambda_\ell, 0) \in \mathcal{S}$. Let $\mathcal{X}_\ell$ denote the maximal connected component of $\mathcal{S}$ containing $(\lambda_\ell, 0)$. Then at least one of the following holds: (i) $\mathcal{X}_\ell$ is unbounded in $\mathbb{R} \times E_{\mathcal{G}} \times \mathbb{R}^2$; (ii) $\mathcal{X}_\ell$ contains $(\lambda_*, 0)$, with $\lambda_* \neq \lambda_\ell$.*

**Proof.** By the linearization principle of Leray and Schauder, (6.21) and (6.22) insure that the Leray Schauder degree of the mapping $v \mapsto v - \mathcal{K}(\lambda, v)$ along the trivial solution, changes sign as $\lambda$ crosses from one side of $\lambda_\ell$ to the other, in some sufficiently small neighborhood of $\lambda_\ell$. The rest of the proof is due to Rabinowitz [35], except that here, because of our definition of $\mathcal{K}(\cdot)$ on $\mathbb{R} \times (\mathcal{B}_M \cap E_{\mathcal{G}}) \times \mathbb{R}^2$, we obtain a third non-exclusive alternative: (iii) $\mathcal{X}_\ell \cap (\mathbb{R} \times \partial\mathcal{B}_M) \neq \varnothing$ for each $M > 0$. But if (iii) is true for all $M > 0$, then (ii) holds. On the other hand, if (iii) does not hold, say, for some sufficiently large $M > 0$, then (i) and/or (ii) remain valid. □

**Remark 6.4.** The crossing condition (6.22) is weaker than (6.12), cf. (6.17). As such, $\mathcal{X}_\ell$ need not contain a smooth local curve of solutions without imposing (6.12). Nonetheless, it's straightforward to show, using Proposition 5.2, that for any sequence of solutions $\{(\lambda_k, v_k)\} \subset \mathcal{X}_\ell$ with $(\lambda_k, v_k) \to (\lambda_\ell, 0)$ in $\mathbb{R} \times E_{\mathcal{G}} \times \mathbb{R}^2$ as $j \to \infty$, there is a subsequence $\left\{ v_{k_n} / \left\| v_{k_n} \right\| \right\} \subset E_{\mathcal{G}} \times \mathbb{R}^2$ converging to $c\hat{z}_\ell$, $c \neq 0$. We also point out that Theorem 6.3 is valid here for any odd-dimensional null space, i.e., the dimension in (6.10) is odd.

There is a well-known strategy for the selection of subgroups in equivariant bifurcation problems, "generically" leading to one-dimensional null spaces [14]. Presuming an appropriate Liapunov-Schmidt reduction, a finite-dimensional bifurcation problem defined on some open neighborhood of, say, $\mathbb{R} \times \mathcal{N} \to \mathcal{N}$, where $\mathcal{N}$ is isomorphic to the null space, is equivariant under the restriction of the group action to $\mathcal{N}$, cf. [5]. In particular, for generic bifurcation problems in the presence of $O(3)$ symmetry, $\mathcal{N} = span\{\rho_{\ell,m}(\mathbf{x}), \ m = -\ell, ..., \ell\}$, $\ell \in \mathbb{N}$, cf. (5.6), and the group action on $\mathcal{N}$, inherited from (6.1), is

$$\sum_{-\ell \le m \le \ell} c_m \rho_{\ell,m}(\mathbf{G}^T \mathbf{x}) \equiv \sum_{-\ell \le m \le \ell} \sum_{j=1}^{2\ell+1} (T_G^\ell)_{mj} c_j \rho_{\ell,m}(\mathbf{x}) \text{ for all } \mathbf{G} \in O(3). \tag{6.23}$$

The representation $\tilde{O}_\ell(3) = \{T_G^\ell : \mathbf{G} \in O(3)\}$ is irreducible, cf. [43]. For each $\ell \in \mathbb{N}$, the maximal isotropy subgroups $\mathcal{G} \subset O(3)$ having 1-dimensional fixed-point spaces, i.e., such that

$$\dim\{\mathbf{c} \in \mathbb{R}^{2\ell+1} : T_G^\ell \mathbf{c} = \mathbf{c} \text{ for all } \mathbf{G} \in \mathcal{G}\} = 1, \tag{6.24}$$

are then enumerated. Fortunately a complete classification of such is provided in [14], cf. also [5]. Bifurcating solutions to the finite-dimensional problem are now readily obtained courtesy of (6.24) ("generically" assuming a version of (6.11) on the 1-dimensional subspace of $\mathcal{N}$). We remark that this strategy, sometimes called the equivariant branching lemma, is not exhaustive in the $O(3)$- symmetry setting. That is, there are subgroups of $O(3)$ for which the dimension in (6.24) is greater than one (for a given $\ell \in \mathbb{N}$) such that the finite-dimensional problem admits generic local solutions characterized by sub-maximal isotropy [6]. We briefly comment on this within the context of global bifurcation in the final section of the paper.

As a starting point here, we read off the appropriate subgroups from [14;Thm 9.9]. However, we still need to address the presence of the 3-dimensional null space $\mathcal{N}_o$, cf. (5.15), (5.17), which is independent of $\lambda$. First, we claim that the irreducible representation for $\ell = 1$ is the same as $O(3)$ itself, viz., $\tilde{O}_1(3) \equiv O(3)$. Indeed by virtue of (5.6), (5.7), observe that

$$x_1 \equiv \rho_{1,1}(\mathbf{x}), x_2 \equiv \rho_{1,-1}(\mathbf{x}), x_3 \equiv \rho_{1,0}(\mathbf{x}), \text{ for all } \mathbf{x} \in S^2. \tag{6.25}$$

Then defining $\mathbf{c} := (c_1^1, c_{-1}^1, c_0^1)$, we may express (5.11) as

$$\upsilon_1 = c_1^1 \rho_{1,1}(\mathbf{x}) + c_{-1}^1 \rho_{1,-1}(\mathbf{x}) + c_0^1 \rho_{1,0}(\mathbf{x}) = c_1^1 x_1 + c_{-1}^1 x_2 + c_0^1 x_3 := \mathbf{c} \cdot \mathbf{x}, \tag{6.26}$$

and thus,

$$\begin{aligned}
c_1^1 \rho_{1,1}(\mathbf{G}^T \mathbf{x}) + c_{-1}^1 \rho_{1,-1}(\mathbf{G}^T \mathbf{x}) + c_0^1 \rho_{1,0}(\mathbf{G}^T \mathbf{x}) &= \mathbf{c} \cdot \mathbf{G}^T \mathbf{x} = \mathbf{Gc} \cdot \mathbf{x} \\
&= \mathbf{Gc} \cdot (\rho_{1,1}(\mathbf{x}), \rho_{1,-1}(\mathbf{x}), \rho_{1,0}(\mathbf{x})) \text{ for all } \mathbf{G} \in O(3).
\end{aligned} \tag{6.27}$$

From (6.23) we conclude that $T_G^1 \equiv \mathbf{G}$ for all $\mathbf{G} \in \mathrm{O}(3)$. This leads to

**Proposition 6.5.** *Given an isotropy subgroup $\mathcal{G} \subset O(3)$, suppose that the fixed-point space in $\mathbb{R}^3$ is trivial, that is*

$$T_G^1 \mathbf{c} \equiv \mathbf{Gc} = \mathbf{c} \text{ for all } \mathbf{G} \in \mathcal{G} \Rightarrow \mathbf{c} = \mathbf{0}. \tag{6.28}$$

*Then the following additions to Proposition 5.1 hold:*

$$\begin{aligned}
\mathcal{N}\left( L(\lambda)|_{X_\mathcal{G} \times \mathbb{R}^2} \right) &= \{0\}, \ \lambda \neq \lambda_\ell; \\
\mathcal{N}\left( L(\lambda_\ell)|_{X_\mathcal{G} \times \mathbb{R}^2} \right) &\subseteq \mathcal{N}_\ell,
\end{aligned} \tag{6.29}$$

**Proof.** Suppose that $(0, \upsilon_1, 0, 0) \in \mathcal{N}_r \cap X_\mathcal{G} \times \mathbb{R}^2$, with $\upsilon_1$ given by (6.26). Then (6.1), (6.7) and (6.27) imply that $\mathbf{Gc} \cdot \mathbf{x} = \mathbf{c} \cdot \mathbf{x}$ for all $\mathbf{G} \in \mathcal{G}, \mathbf{x} \in S^2$, i.e., the left side of (6.29) holds. Thus, $\mathbf{c} = (c_1^1, c_{-1}^1, c_0^1) = \mathbf{0} \Rightarrow \upsilon_1 \equiv 0$. □

Condition (6.28) requires that there are no invariant lines through the origin in $\mathbb{R}^3$ under the group action of $\mathcal{G}$. This is indeed the case for all but one class of problems associated with the isotropy subgroup classification of [14], namely those having axisymmetric solutions corresponding to odd $\ell \in \mathbb{N}$. We discuss a simple resolution for obtaining the existence of those solutions as well at the end of this section, cf. Remark 6.7. With this in hand, we now summarize:

**Proposition 6.6.** *If $\lambda_\ell$ is a root of the characteristic equation (5.9) for a given $\ell \in \mathbb{N}$, and if (6.22) is satisfied, then the global bifurcation Theorem 6.3 is valid for each of the maximal isotropy subgroups listed in [14; Thm.9.9]. In addition, if (6.11) holds, then the global solution continuum $\mathcal{X}_\ell$ contains the local bifurcation curve (6.13).*

We now illustrate with several concrete examples. In what follows, we employ the subgroup notation of [14]: The *planar* subgroups, denoted $SO(2), O(2), D_n, Z_n \, (n = 2, 3, ...) \subset SO(3)$, are oriented in the obvious way with respect to the $x_3$-axis. We specify the orientation of the non-planar or exceptional subgroups, $\mathbb{T}, \mathbb{O}, \mathbb{I} \subset SO(3)$, as in [15]. Specifically, we choose

$$\boldsymbol{\Gamma}_1 = \begin{pmatrix} 0 & 1 & 0 \\ 0 & 0 & 1 \\ 1 & 0 & 0 \end{pmatrix} \text{ and } \boldsymbol{\Gamma}_2 = \begin{pmatrix} -1 & 0 & 0 \\ 0 & -1 & 0 \\ 0 & 0 & 1 \end{pmatrix} \tag{6.30}$$

as the generators of the *tetrahedral* group $\mathbb{T}$. We then choose a copy of the *octahedral* group $\mathbb{O}$ so that it is a supergroup of $\mathbb{T}$, as specified above. The *icosahedral* group $\mathbb{I}$ is chosen so that it contains two $\mathbb{Z}_5$ subgroups – one about the $(0,0,1)$ axis and one about $(-2/\sqrt{5}, 0, 1/\sqrt{5})$. We also observe that $O(3) = SO(3) \oplus \mathbb{Z}_2^c$, where $\mathbb{Z}_2^c := \{I, -I\}$.

One last issue needs to be addressed before presenting our examples. Unlike the calculations leading to (6.24) (which involve traces only, cf. [5], [14]), the full representation $\tilde{O}_\ell(3) = \{T_G^\ell : \mathbf{G} \in O(3)\}$ (its generators) is needed in the actual construction of the *bifurcation direction,*

$$\hat{z}_\ell = (\hat{\rho}_\ell, \tau_\ell \hat{\rho}_\ell, 0, 0), \tag{6.31}$$

for a given $\ell \in \mathbb{N}$ and isotropy subgroup $\mathcal{G} \subset O(3)$. Specifically, we seek $\hat{\mathbf{c}} \in \mathbb{R}^{2\ell+1}$ such that

$$T_G^\ell \hat{\mathbf{c}} = \hat{\mathbf{c}} \text{ for all } \mathbf{G} \in \mathcal{G}. \tag{6.32}$$

which yields

$$\hat{\rho}_\ell = \sum_{-l \le m \le l} \hat{c}_m \rho_{\ell,m}. \tag{6.33}$$

While the "directions" (6.33) have not been exhaustively catalogued, some have been determined and presented in our earlier papers [15], [19], of which we take full advantage here. We present details for some that seem to be quite relevant to the experimentally obtained patterns presented in [4]. We mention that perhaps a more direct approach than (6.32) comes from the classical method in [34], restricted here to the null spaces $\mathcal{N}_\ell$. We refer to [9] for details, which will also appear elsewhere.

**Examples:** We point out that the first example – in particular (6.34) – corrects an error associated with that presented in [19].

I. $\ell = 3, \mathcal{G} = D_6^d$ :

Here $D_6^d \not\subset SO(3)$ denotes the group generated by

$$\mathbf{\Gamma}_1 = \begin{pmatrix} \cos(2\pi/3) & \sin(2\pi/3) & 0 \\ -\sin(2\pi/3) & \cos(2\pi/3) & 0 \\ 0 & 0 & 1 \end{pmatrix}, \mathbf{\Gamma}_2 = \begin{pmatrix} 1 & 0 & 0 \\ 0 & -1 & 0 \\ 0 & 0 & -1 \end{pmatrix} \text{ and } \mathbf{\Gamma}_2 = \begin{pmatrix} 1 & 0 & 0 \\ 0 & 1 & 0 \\ 0 & 0 & -1 \end{pmatrix}. \tag{6.34}$$

$\mathcal{N}\left(L(\lambda_3)\big|_{X_G \times \mathbb{R}^2}\right) = span\{\hat{z}_3 = (\rho_{3,3}, \tau_3 \rho_{3,3}, 0, 0)\}$. The nodal set of $\rho_{3,3}$, cf. (5.6), is depicted below in Figure 6.1(a).

II. $l = 3, \mathcal{G} = \mathbb{O}^-$ :

Here $\mathbb{O}^- \not\subset SO(3)$ denotes the subgroup $\mathbb{O}^- = \mathbb{T} \cup \{-\mathbf{\Gamma} : \mathbf{\Gamma} \in \mathbb{O} \setminus \mathbb{T}\}$, for the copy of $\mathbb{T}$ given above. $\mathcal{N}\left(L(\lambda_3)\big|_{X_G \times \mathbb{R}^2}\right) = span\{\hat{z}_3 = (\rho_{3,-2}, \tau_3 \rho_{3,-2}, 0, 0)\}$. The nodal set of $\rho_{3,-2}$ is depicted below in Figure 6.1(b).

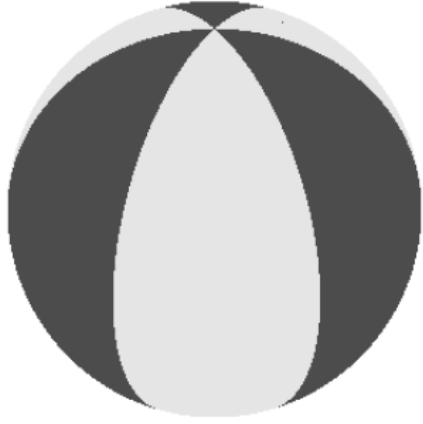

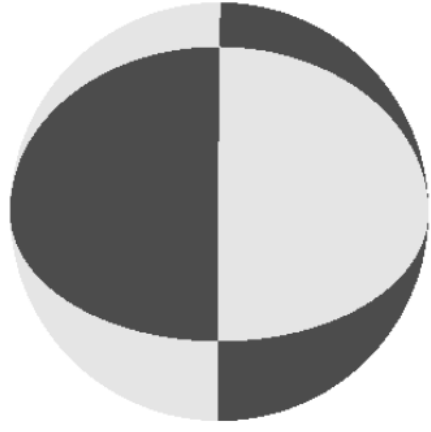

(a) $\mathcal{G} = D_6^d$.  (b) $\mathcal{G} = \mathbb{O}^-$.

Figure 6.1. Nodal sets of eigenfunctions corresponding to bifurcations for $\ell = 3$.

III. $\ell = 4, \mathcal{G} = O(2) \oplus \mathbb{Z}_2^c$ :

$\mathcal{N}\left(L(\lambda_4)|_{X_\mathcal{G} \times \mathbb{R}^2}\right) = span\{\hat{z}_4 = (\rho_{4,0}, \tau_4 \rho_{4,0}, 0, 0)\}$. The nodal set of $\rho_{4,0}$ is depicted below in Figure 6.2(a).

IV. $\ell = 4, \mathcal{G} = \mathbb{O} \oplus \mathbb{Z}_2^c$ :

$\mathcal{N}\left(L(\lambda_4^\beta)|_{X_\mathcal{G} \times \mathbb{R}^2}\right) = span\{\hat{z}_4 = (\hat{\rho}_4, \tau_4 \hat{\rho}_4, 0, 0)\}, \hat{\rho}_4 = 168\rho_{4,0} + \rho_{4,4}$. The nodal set of $\hat{\rho}_4$ is depicted below in Figure 6.2(b).

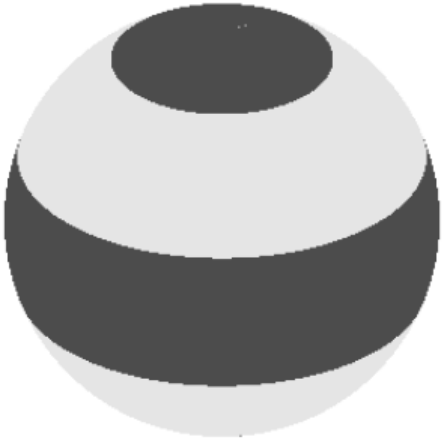

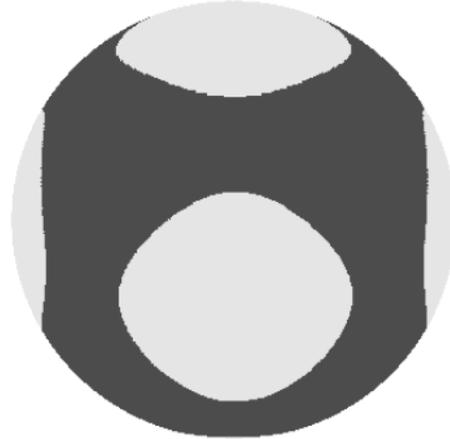

(a) $\mathcal{G} = O(2) \oplus \mathbb{Z}_2^c$.  (b) $\mathcal{G} = \mathbb{O} \oplus \mathbb{Z}_2^c$.

Figure 6.2. Nodal sets of eigenfunctions corresponding to bifurcations for $\ell = 4$.

**V.** $\ell = 6, \mathcal{G} = \mathbb{I} \oplus \mathbb{Z}_2^c :$

$\mathcal{N}\left( L(\lambda_6)|_{X_{\mathcal{G}} \times \mathbb{R}^2} \right) = span\{\hat{z}_6 = (\hat{\rho}_6, \tau_6 \hat{\rho}_6, 0, 0)\}, \hat{\rho}_6 = 3,960\rho_{6,0} - \rho_{6,5}.$ The nodal set of $\hat{\rho}_6$ is depicted below in Figure 6.3(a).

**VI.** $\ell = 10, \mathcal{G} = \mathbb{I} \oplus \mathbb{Z}_2^c :$

$\mathcal{N}\left( L(\lambda_{10})|_{X_{\mathcal{G}} \times \mathbb{R}^2} \right) = span\{\hat{z}_{10} = (\hat{\rho}_{10}, \tau_{10} \hat{\rho}_{10}, 0, 0)\}, \hat{\rho}_{10} = \rho_{10,10} + 896,313,600\rho_{10,0} + 27,360\rho_{10,5}.$ The nodal set of $\hat{\rho}_{10}$ is depicted below in Figure 6.3(b).

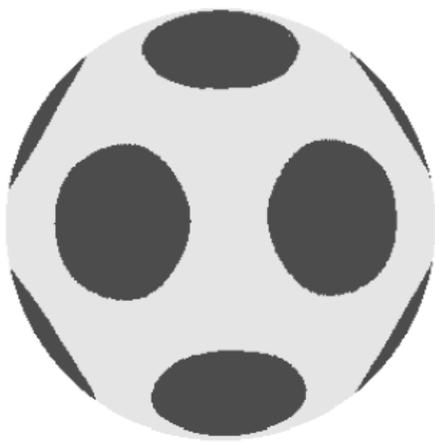 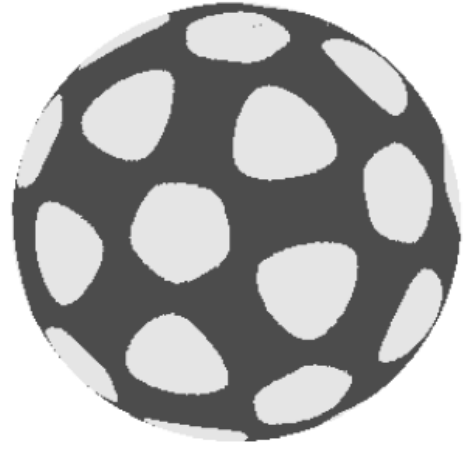

(a) $\mathcal{G} = \mathbb{I} \oplus \mathbb{Z}_2^c$ for $\ell = 6$.      (b) $\mathcal{G} = \mathbb{I} \oplus \mathbb{Z}_2^c$ for $\ell = 10$.

Figure 6.3. Nodal sets of eigenfunctions corresponding to bifurcations.

**VIII.** $\ell = 12, \mathcal{G} = \mathbb{I} \oplus \mathbb{Z}_2^c :$

$\mathcal{N}\left( L(\lambda_{12})|_{X_{\mathcal{G}} \times \mathbb{R}^2} \right) = span\{\hat{z}_{12} = (\hat{\rho}_{12}, \tau_{12} \hat{\rho}_{12}, 0, 0)\}, \hat{\rho}_{12} = 57,001,190,400\rho_{12,0} - 221,760\rho_{12,5} + 4\rho_{12,10}.$ The nodal set of $\hat{\rho}_{12}$ is depicted below in Figure 6.4.

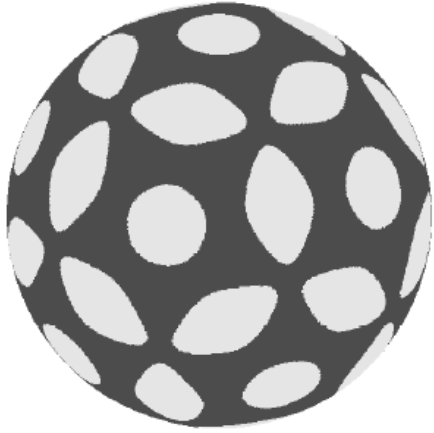

$$\mathcal{G} = \mathbb{I} \oplus \mathbb{Z}_2^c$$

Figure 6.4. Nodal set of the eigenfunction corresponding to bifurcation for $\ell = 12$.

**Remark 6.7** The pertinent isotropy subgroup for axisymmetric solutions corresponding to odd $\ell \in \mathbb{N}$ is $O(2)^-$, which is generated by $SO(2)$ and a reflection across the $x_2$-$x_3$ plane. Obviously (6.28) does not hold in this case; the entire $x_3$ axis is fixed by the group action. This, in turn, implies that $\upsilon(\mathbf{x}) = c_0^1 \rho_{1,0}(\mathbf{x}) \equiv c_0^1 \cos\theta$, cf. (5.7), representing an arbitrary vertical translation, is present in the null space of the linear operator. Of course in this case, $v \in X_{O^-(2)} \Rightarrow (\phi(\mathbf{x}), u(\mathbf{x})) = (\tilde{\phi}(\theta), \tilde{u}(\theta))$, and the null solution is readily factored out by appending, e.g., the boundary condition $\tilde{u}(0) = 0$ to the governing 2-point boundary value problem.

## 7. Concluding Remarks

The basic necessary condition for bifurcation, viz., that the characteristic equation (5.9) has a root $\lambda_\ell$ for a given $\ell \in \mathbb{N}$, involves the second derivative of the constitutive function $\Psi = W + B + E$, cf. (2.1) and (3.37)₃. In the special case that the bending moduli, $B$ and $E$, are constants, this obviously involves only the phase-field potential $W$, and there are precisely two such roots for each $\ell \in \mathbb{N}$, for $\varepsilon$ sufficiently small, each satisfying (6.12). This follows directly from (5.9) and the assumed graph of $W''$ depicted in Figure 2.1. In the more general case considered here, it is reasonable to assume that both $B(\cdot)$ and $E(\cdot)$ are essentially constant outside of the spinodal region, with a smooth transition within that interval $(m_1, m_2)$. As such, the function $\Psi''(\cdot)$ takes on negative values in the spinodal region only, and again, for sufficiently small $\varepsilon$, (5.9) has at least one root for each $\ell \in \mathbb{N}$. Thus, as $\varepsilon$ approaches zero from above, and presuming that either (6.12) or (6.22) holds, then Proposition 6.6 is valid, i.e., we have the existence of a global branch of bifurcating solutions correspond to *each* of the maximal isotropy subgroups of $O(3)$ as classified in [14; Thm 9.9].

As mentioned in Section 6, certain local, generic, $O(3)$- symmetry-breaking solutions characterized by sub-maximal isotropy have been uncovered [6]. These follow from an analysis of the reduced bifurcation equations associated with sub-maximal isotropy subgroups. In particular, this engenders higher-dimensional fixed-point eigen-subspaces, i.e., the right side of (6.10) is greater than one, which corresponds to the number of reduced bifurcation equations. As in [19], the results of [6] suggest the existence of *local* solutions to our problem (although we have not pursued the details here). In cases when the null space is odd-dimensional, our global bifurcation theorem is valid, cf. Remark 6.4. However, this could give nothing new without verifying the local analyses in [6]; as pointed out in [19], every such fixed-point space already contains one of the global solutions branches given by Proposition 6.6. Moreover, even if the local results of [6] are valid here, it does *not* then follow that our problem possesses global solutions (outside of a small neighborhood of a given bifurcation point) that are precisely characterized by the associated sub-maximal isotropy. We refer to [19] for a detailed discussion.

Our local analysis, staring from (5.15)-(5.17) and culminating in (6.14) or, more generally, Remark 6.4, illuminates the role of the bending moduli functions, $B(\cdot)$ and $E(\cdot)$, in yielding a phenomenon not seen in the constant case. Specifically, suppose that $\hat{z}_\ell$ is the local bifurcation direction according to (6.31)-(6.33). Then by (6.14)$_3$, (6.31) and (6.33), we see that the displacement or shape variable $\upsilon_\ell = \tau_\ell \hat{\rho}_\ell$ is *not* present to first order in (6.14) when the bending moduli are constant, which is due to the fact that $\tau_\ell \equiv 0, \ell \geq 2$, cf. (5.14), (5.15). Recall that $\tau_1 = 0$, regardless of the behavior of the bending moduli. Indeed, as addressed in Proposition 6.5 and Remark 6.7, $\upsilon_1 \neq 0$ represents a rigid translation. In any case, from (5.14), we see that the non-vanishing of either $B'(\lambda_\ell)$ and/or $E'(\lambda_\ell)$ within the spinodal region $(m_1, m_2)$, which is consistent with the discussion in the paragraph above, typically engenders $\tau_\ell \neq 0, \ell \geq 2$. In general for $\ell = 1$, or when both $B$ and $E$ are constants for all $\ell = 2, 3, \ldots$, the local bifurcating branch (6.14) is the same - to first order - as that for the phase field on a fixed sphere, as treated in [19]. This raises the question as to the possible existence of nontrivial solutions characterized by a non-constant phase field $\phi$ on the undeformed sphere here in this problem, i.e., with $u \equiv 0$, cf. (3.1). We claim this is not possible. Indeed, return to (2.12) with $\mathbf{L} = -\mathbf{I}_x$, cf. (3.8), corresponding to $H \equiv -1$ and $K \equiv 1$; we find that $W(\phi) - \mu\phi = \gamma + p/2$, i.e. $\phi$ is a constant. But then by virtue of (2.2)$_2$, $\phi \equiv \lambda$, and we are back to the trivial solution (3.36), (3.37). We conclude that the uniform (trivial) solution $\equiv$ the only possible (smooth) solution of our problem on the undeformed sphere.

Our approach in this work can also be used to establish existence of global, non-axisymmetric bifurcating solutions for the classical Helfrich model. The latter corresponds to (2.1) in the absence of the phase field, engendering constant bending moduli, and subject to (2.2)$_1$; the second term in the integrand in (2.1) now falls out, by virtue of the Gauss-Bonnet theorem. The pressure is the natural bifurcation parameter, in which case we know from the literature that axisymmetric buckled states exist – these bifurcations occur at negative values of the internal pressure - corresponding to external compressive pressure, cf. [25], [32]. At any bifurcation point (where an axisymmetric solution always exists), our analysis here provides global solution branches characterized by spatial symmetries as in Sections 6. In the same way we can analyze global "buckling" solutions here in our problem by fixing the phase-field parameter $\lambda$ and treating "$- p$" as the bifurcation parameter. Similarly, we can fix both $p$ and $\lambda$ and obtain global bifurcation results in the reciprocal of the small parameter, viz., $1/\varepsilon$, cf. (3.33), as carried out in [20].

Similar to the results in [19], most, if not all *local* bifurcating solutions are unstable, and our proof of global existence, based upon the Leray-Schauder degree, provides no information about the stability of global solutions. Nonetheless, we expect that many of our global solution branches are totally bounded and connected (in the Cartesian product of the parameter line and the Banach solution space), as depicted schematically in Figure 6 of [19], strongly suggesting stability "far" from the trivial solution. Of course, the experimental results of [4] and the numerical, gradient flow results of [12] and [37] also suggest the stability of solutions with symmetries like those depicted in Figures 6.3, 6.4, among others. Stability can also be checked on global solutions via numerical computation of the second variation, e.g., as in [22].

Precisely the same issues discussed at the start of Section 3, viz., in-plane fluidity and the degeneracy associated with a Lagrangian description, also plague numerical computation in models of lipid-bilayer vesicles, cf. [12], [13], [30], [38]. Typically some type of added in-plane stiffness for static simulations [13] and/or in-plane dissipation for gradient-flow dynamics [12], [30], [38] is required. In this regard, a radial-map description like (3.1) is well suited to numerical computation of equilibria; the grossly underdetermined in-plane deformation is absent in the resulting formulation.

Our results presented herewith are perhaps limited by the radial-graph formulation. That is, there may exist more general equilibria of (2.2), (2.11), (2.12) than those captured by (3.1). However, it is worth noting that all of the observed (apparent) equilibria reported in [4] appear as radial graphs over the sphere.

## Acknowledgements

This work was supported in part by the National Science Foundation through grants DMS-1007830 and DMS-1312377, which is gratefully acknowledged. We also thank the following people for valuable comments, criticisms and suggestions throughout the long gestation period of this work: Tobias Baumgart, Sovan Das, Luca Deseri, Xiaodong Cao, Jim Jenkins, Alexander Mielke, Phoebus Rosakis, Siming Zhao and Giuseppe Zurlo.